\documentclass{article}

\usepackage{amsmath}
\usepackage{amsfonts}
\usepackage{amscd}
\usepackage{amssymb}

\input xypic

\newtheorem{defn}{Definition}[section]
\newtheorem{thm}[defn]{Theorem}
\newtheorem{prop}[defn]{Proposition}
\newtheorem{lemma}[defn]{Lemma}

\newtheorem{claim}[defn]{Claim}
\newtheorem{cor}[defn]{Corollary}

\newtheorem{rem}[defn]{Remark}

\newcommand{\lm}{\ensuremath{\longrightarrow}}

\newcommand{\Gm}{\ensuremath{\mathbb{G}_{\text{m}}}}
\DeclareMathOperator{\Obj}{\mbox{Obj}\,}

\DeclareMathOperator{\Hom}{\mbox{Hom}}

\DeclareMathOperator{\shom}{\ensuremath{\mathcal{H}\mathit{om}}}

\DeclareMathOperator{\send}{\ensuremath{\mathcal{E}\!\mathit{nd}}}
\DeclareMathOperator{\End}{\mbox{End}}

\DeclareMathOperator{\Ext}{\mbox{Ext}}

\DeclareMathOperator{\id}{\mbox{id}\,}
\DeclareMathOperator{\cd}{\mbox{cd}\,}
\DeclareMathOperator{\im}{\mbox{im}\,}
\DeclareMathOperator{\spec}{\mbox{Spec}\,}
\DeclareMathOperator{\proj}{\mbox{Proj}\,}
\DeclareMathOperator{\mo}{\mbox{mod}}
\DeclareMathOperator{\Mod}{\mbox{Mod}}
\DeclareMathOperator{\gr}{\mbox{gr}}

\DeclareMathOperator{\Gr}{\mbox{Gr}}
\DeclareMathOperator{\Br}{\mbox{Br}\,}
\DeclareMathOperator{\Pic}{\mbox{Pic}\,}

\DeclareMathOperator{\coker}{\mbox{coker}\,}

\DeclareMathOperator{\Hilb}{\mbox{Hilb}}
\DeclareMathOperator{\Aff}{\mbox{NAff}}

\DeclareMathOperator{\qcoh}{\mbox{QCoh}}
\DeclareMathOperator{\mmm}{\mathfrak{m}}
\DeclareMathOperator{\s}{\sigma}

\DeclareMathOperator{\w}{\omega}

\DeclareMathOperator{\PP}{\mathbb{P}}
\DeclareMathOperator{\bu}{\bullet}
\DeclareMathOperator{\Z}{\mathbb{Z}}
\DeclareMathOperator{\calm}{\mathcal{M}}
\DeclareMathOperator{\calo}{\mathcal{O}}
\DeclareMathOperator{\caln}{\mathcal{N}}
\DeclareMathOperator{\call}{\mathcal{L}}

\DeclareMathOperator{\cala}{\mathcal{A}}
\DeclareMathOperator{\calb}{\mathcal{B}}
\DeclareMathOperator{\calc}{\mathcal{C}}

\DeclareMathOperator{\caly}{\mathcal{Y}}

\DeclareMathOperator{\oy}{\mathcal{O}_{Y}}

\DeclareMathOperator{\Tt}{\ensuremath{\tilde{T}}}
\DeclareMathOperator{\Ht}{\ensuremath{\tilde{H}}}

\begin{document}

\begin{center}
\LARGE \textbf{Twisted rings and moduli stacks of ``fat'' point modules in non-commutative projective geometry}
\end{center}

\begin{center}
  DANIEL CHAN
\footnote{This project was supported by the Australian Research Council, Discovery Project Grant DP0880143.}
\end{center}

\begin{center}
  {\em University of New South Wales}
\end{center}

\begin{center}
e-mail address:{\em danielc@unsw.edu.au}
\end{center}

\begin{abstract}
The Hilbert scheme of point modules was introduced by Artin-Tate-Van den Bergh to study non-commutative graded algebras. The key tool is the construction of a map from the algebra to a twisted ring on this Hilbert scheme. In this paper, we study moduli stacks of more general ``fat'' point modules, and show that there is a similar map to a twisted ring associated to the stack. This is used to provide a sufficient criterion for a non-commutative projective surface to be birationally PI. It is hoped that such a criterion will be useful in understanding Mike Artin's conjecture on the birational classification of non-commutative surfaces.
\end{abstract}

Throughout, all objects and maps are assumed to be defined over some algebraically closed base field $k$. Also, the unadorned tensor product symbol $\otimes$ will be used for tensor products over $k$. 

\vspace{2mm}
\section{Introduction}

A natural way to study a non-commutative algebra is to analyse its moduli of simple modules. In non-commutative projective geometry, one studies a non-commutative connected graded algebra $A$ and we note of course, the simple objects in the category $A-\Gr$ of graded $A$-modules give no information about $A$. One of the insights of Artin, Tate and Van den Bergh, is to consider instead the category of $A$-tails $\proj A := A-\Gr /\text{tors}$ where tors is the subcategory of $A_{>0}$-torsion modules. In [ATV], they study the Hilbert scheme of certain simple objects in $\proj A$ called point modules. In this paper, we generalise their methods by introducing Artin stacks to study the moduli of more general simple objects sometimes dubbed ``fat'' points.


To be more precise, we start with a noetherian graded $k$-algebra, finitely generated in degree one which satisfies good homological properties (see subsection~\ref{ssextraass}) the most important of which is strong $\chi$. This ensures the existence of Hilbert schemes \`a la Artin-Zhang [AZ01]. For an integer $m\geq 1$, an $m$-point is an $A$-tail $M$ such that the corresponding (saturated) $A$-module has constant Hilbert function $m$, that is, in the language of the non-commutative projective geometry of Artin-Zhang [AZ94] (see section~\ref{smodulifat}), $\dim_k H^0(M(i)) = m$ for all $i \in \Z$. The point modules of [ATV] are essentially thus the 1-points, while the $m>1$ case corresponds to what Artin has called a ``fat'' point. We show that $m$-points are ``parametrised'' by an Artin stack. Furthermore, the simple $m$-points form an open substack $\caly$. As usual, we construct this by studying an appropriate Hilbert scheme. We show that there is an open subset $H$ of a projective Hilbert scheme with a free action of $PGL_m$ on it, such that $\caly$ is isomorphic to the quotient stack $[H/GL_m]$ where we have extended the $PGL_m$-action trivially to $GL_m$. For point modules, the stack is just $[H/\mathbb{G}_m]$ where $\mathbb{G}_m$ acts trivially so one can get away with working only with the Hilbert scheme as is done in [ATV], [RZ] and elsewhere. One of the points of this paper, is to make the observation that only by considering the stack of $m$-points can one generalise the methods of [ATV]. 

We now recall the utility of twisted rings in [ATV], later generalised by Rogalski and Zhang in [RZ]. 
There is a shift functor on $\proj A$ induced by shifting graded modules, namely, given $M \in A-\Gr$ we define $M(1) \in A-\Gr$ to be the graded module with $i$-th graded piece $M(1)_i := M_{i+1}$.  
Let $Y$ be the Hilbert scheme of point modules, or point scheme for short, and $\calm$ be the universal point module on $Y$. The universal property ensures that the shift functor induces an automorphism $\s$ of $Y$, such that $\calm(1) \simeq \s^*\calm \otimes_Y \call$ for some line bundle $\call$ on $Y$. We can form a twisted sheaf of rings 
$$\mathcal{E}:=\oy \oplus \call \oplus (\s^* \call \otimes_Y \call) \oplus \ldots .$$
Alternatively, one can write the $d$-th graded component of this ring as $\call_{\s}^{\otimes d}$ where $\call_{\s}$ is an invertible $\oy$-bimodule as defined in [AV]. The twisted ring above is very useful for studying $A$ since there is a canonical map from $A \lm H^0(\mathcal{E})$. In fact, $\proj H^0(\mathcal{E})$ is essentially the ``maximal commutative subscheme'' of $\proj A$. We note that this does not necessarily mean $H^0(\mathcal{E})$ is $A$ modulo the commutator ideal since $\proj A$ can be commutative even if $A$ is ``highly'' non-commutative. 

We seek a similar construction for $m$-points. Again, there is a universal $m$-point, say $\calm$ on $\caly$ and an automorphism $\s: \caly \lm \caly$ induced by the shift functor such that $\calm(1) \simeq \s^*\calm$. 
We associate to these data, a twisted ring by passing to an Azumaya algebra as follows. 
First note that the coarse moduli space of $\caly$, is the separated algebraic space $Y = H/PGL_m$. The $PGL_m$-torsor $p:H \lm Y$ corresponds to an Azumaya algebra $\cala$ on $Y$ which is our $m$-point analogue of $\oy$. The automorphism $\s$ corresponds to an invertible $\cala$-bimodule $\calb$ which takes the role of $\call_{\s}$ in the point scheme case. We can thus form the twisted sheaf of rings $\mathcal{E}$ as before and furthermore, there is a map $A \lm H^0(\mathcal{E})$. We like to think of the twisted ring $H^0(\mathcal{E})$ as capturing the stratum of $\proj A$ which is PI of PI-degree exactly $m$. In particular, if $\proj A$ arises from an order on a variety, then one expects that for appropriate $m$, $Y$ is the Azumaya locus of the order and $\cala$ is the restriction of that order to $Y$. 

Our interest in studying $m$-points  stems from Artin's conjecture [A97] on the classification of non-commutative surfaces which we paraphrase imprecisely as follows. Suppose $A$ has Gelfand-Kirillov dimension three so $\proj A$ is a non-commutative surface. Then Artin conjectures that the degree zero part $Q(A)_0$ of the graded quotient ring $Q(A)$ is either i) (birationally ruled) isomorphic to $Q(A')_0$ where $A'$ gives a non-commutative ruled surface as defined by Van den Bergh and Patrick [VdB01p], [Pat] or, ii) (birationally PI) finite over its centre. 
Our main result is theorem~\ref{tbiratPI} which roughly states that a non-commutative surface is birationally PI if it has a surface worth of simple $m$-points. This reinforces the well-observed phenomenon that ``points disappear the more non-commutative one gets''. The proof of the theorem is based on the map to the twisted ring above and the fact that the twisted ring is birationally PI.

We hope this paper is readable by, and of interest to both algebraists and algebraic geometers. Consequently, we have included more expository material than is usual. Section~\ref{sazu} contains some background material on the particular stacks we are interested in. In section~\ref{sendos}, we look at $A$-modules $\calm_{\bullet}$ with an isomorphism $z:\calm_{\bullet}(1) \stackrel{\sim}{\lm} \s^* \calm_{\bullet}$ like the universal simple $m$-point mentioned above, and identify the endomorphisms compatible with the shift isomorphism $z$ as a twisted ring. In section~\ref{smodulifat}, we examine the moduli stack of $m$-points and more generally, of 0-dimensional $A$-tails. We prove finally in section~\ref{sbirPI}, our main result, our criterion for being birationally PI. In section~\ref{sloccrit}, we give a local criterion for when a non-commutative projective surface has a surface worth of simple $m$-points. Presumably, one can approach this question from deformation theory too. 

\textbf{Acknowledgements:} I would like to thank James Zhang for sharing his insight during several interesting discussions.

\vspace{2mm}

\section{Stacks associated to Azumaya algebras}  \label{sazu}  

In this section, we will recall a well-known relationship between Azumaya algebras and certain stacks. It should serve to fix notation and as a quick crash course in the relevant theory of stacks for the algebraist. The theory also motivates much of this paper. We refer the reader unfamiliar with algebraic spaces to [A73] and comment only that algebraic spaces are a mild generalisation of schemes, so replacing spaces with schemes on a first reading is but a ``white'' lie. Good references for stacks include [Vis], [LMB] and [Gom]. For the algebraist who might be put off by stacks, we remark that the only stacks used will be those described in this section, the study of which reduces to $GL_m$-equivariant geometry.

Let $\cala$ be an Azumaya algebra of rank $m^2$ over a noetherian separated algebraic space $Y$. If $\cala$ is trivial in the Brauer group $\Br Y$ so say $\cala \simeq \send_Y \calm$, then the simple modules are parametrised by the space $Y$ and $\calm$ is the universal simple. If $\cala$ is non-trivial in the Brauer group, then although the closed points of $Y$ parametrise the simples set-theoretically, there is no universal simple on $Y$. Algebraic geometers say that the simple modules are parametrised by a stack $\caly$ (which we elucidate below) and that this stack has a coarse moduli space $Y$. 

There are two ways of viewing a stack $\caly$. Let $\Aff$ denote the category of noetherian affine schemes. The first viewpoint is as a category fibred in groupoids $\caly \lm \Aff$ which satisfies some extra axioms (see [Vis]). If one picks a cleavage for this stack, then it can be viewed as a pseudo-functor or lax 2-functor $\caly:\Aff \lm \mathcal{C}at$, a point of view which is closer to Grothendieck's functor of points. The data for such a pseudo-functor includes, for each test scheme $T \in \Aff$, a category $\caly(T)$ (corresponding to ``$T$-points'' of the stack) and for each morphism $u:S \lm T$ in $\Aff$, a pull-back ``functor'' $u^*:\caly(T) \lm \caly(S)$. It can fail to be a functor since we only require $\id^*$ to be naturally equivalent to $\id$ and $(uv)^*$ to be naturally equivalent to $v^*u^*$. The reader should consult [Vis] for the other axioms. Stacks are isomorphic if, as fibred categories, they are naturally equivalent. 

The stack $\caly$ of simple left $\cala$-modules is defined by setting $\caly(T)$ to be the category of pairs $(t:T \lm Y,\calm)$ where $\calm$ is a $t^*\cala$-module which is locally free over $T$ of rank $m$. Morphisms are isomorphisms of $t^*\cala$-modules and the pull-back functor is the obvious one.

To study this stack of simples, we utilise the Hilbert scheme $\Hilb (\cala^m,m)$ of quotients of $\cala^m$ which are $\cala$-modules of dimension $m$ over $k$ (and in particular are supported in a zero dimensional set). We consider the open subset $H \subseteq \Hilb (\cala^m,m)$ consisting of the $k$-points $q:\cala^m \lm M$ such that the composite
\[  k^m \lm H^0(\cala^m) \lm H^0(M)  \]
is an isomorphism and $M$ is simple. Now $GL_m$ acts on $\cala^m$ and so on the Hilbert scheme and hence on the invariant subset $H$. Moreover, Schur's lemma ensures that $\Gm \subseteq GL_m$ is the stabiliser of every point. The points of $H$ are given by a simple module $M$ and an identification of $k^m$ with $H^0(M)$ up to scalar. Hence the simples correspond to the $PGL_m$-orbits of $H$ and, since the action is free, there is a GIT quotient $p:H \lm H/PGL_m \simeq Y$.

An alternative viewpoint of the stack of simples is via the stack $\caly'= [H/GL_m]$ which is defined as follows. The category $\caly'(T)$ has objects diagrams 
$$ \diagram  \Tt \rto^u \dto^v & H \\ T & \enddiagram $$
where $v:\Tt \lm T$ is a $GL_m$-torsor and $u:\Tt \lm H$ is a $GL_m$-equivariant map. The reader can easily determine the morphisms in the category $\caly'(T)$ as well as the pull-back functor. 

To see how $\caly$ and $\caly'$ determine the same stack, consider the universal simple quotient 
\[  p^* \cala^m \lm \calm    .\]
Technically, $\calm$ is a sheaf on $H \times Y$, but since it is supported on the graph of $p:H \lm Y$, we will view it as a sheaf on $H$. Suppose first that we have an object of $\caly(T)$ given by a map $t:T \lm Y$ and a $t^*\cala$-module $\caln$. Now $\caln$ is locally free of rank $m$ on $T$ so its frame bundle $\Tt \lm T$ is a $GL_m$-torsor. Also, the $k$-points of $\Tt$ are given by simple modules $M$ and an isomorphism of $k^m \simeq H^0(M)$ and hence, a quotient map $\cala^m \lm M$. The universal property of Hilbert schemes thus furnishes us with a $GL_m$-equivariant map $\Tt \lm H$ and thus an object of $\caly'(T)$. It is easy to see that we thus obtain a functor $\caly \lm \caly'$. 

For the inverse equivalence, consider an object of $\caly'(T)$ given by $v:\Tt \lm T,u:\Tt \lm H$ as above. There is an induced map $t:T \lm Y$ to the coarse moduli space. We need to use Galois descent as follows. The pull-back functor $v^*: \qcoh T \lm \qcoh \Tt$ induces an equivalence of categories between quasi-coherent $T$-modules and quasi-coherent $\Tt$-modules $\caln$ equipped with a $GL_m$-action $\Phi:p_2^* \caln \stackrel{\sim}{\lm} \alpha^* \caln$ where 
$$\diagram  GL_m \times \Tt \rto<1ex>^(.65){p_2} \rto<-1ex>_(.65)\alpha & \Tt
\enddiagram$$ 
are the projection and action maps. Now $T$ is affine so $\Tt$ is too, and the $GL_m$-action is just an action of $GL_m$ on the $H^0(\calo_{\Tt})$-module $H^0(\Tt,\caln)$ which is compatible with the action of $GL_m$ on scalars in $H^0(\calo_{\Tt})$. 
In this affine case, given a quasi-coherent $\Tt$-module $\caln$ with a $GL_m$-action, the corresponding $\calo_T$-module can be recovered as the invariants $\caln^{GL_m}$. Consequently, we will in future write $\caln^{GL_m}$ for the descended module even in the case where $T$ is not affine. Note that $\cala^m$ has a natural $GL_m$-action which induces a $GL_m$-action of the universal simple quotient $\calm$. Moreover, the action is of $p^*\cala$-modules in the sense that $\Phi: p_2^* \calm \stackrel{\sim}{\lm} \alpha^* \calm$ is an isomorphism of $p^*\cala$-modules too. Hence, $u^*\calm$ is a $GL_m$-equivariant $\cala$-module which descends to a $t^* \cala$-module on $T$. We thus obtain a functor $\caly' \lm \caly$. 

A quasi-coherent sheaf $\calm$ on $\caly$ consists of the following data: for each object $\tau \in \caly(T),\ T \in \Aff$ we give a quasi-coherent sheaf $\tau^* \calm \in \qcoh(T)$ which is required to respect pullback. We let $\qcoh(\caly)$ denote the category of quasi-coherent sheaves on $\caly$. Note that the natural quotient map $\pi: H \lm \caly$ gives rise to an object $\pi \in \caly(H)$ so $\pi^*\calm\in \qcoh(H)$. Let $p_2,\alpha: GL_m \times H \lm H$ be the projection and action maps. Then by descent theory, the only data required to determine $\calm \in \qcoh(\caly)$ is $\pi^*\calm$ and the isomorphism $p_2^*\pi^*\calm \simeq \alpha^*\pi^* \calm$, in other words, a $GL_m$-action on $\pi^*\calm$. Thus $\qcoh(\caly)$ is naturally equivalent to the category of $GL_m$-equivariant quasi-coherent sheaves on $H$ and we will often identify the two categories this way. Given a $GL_m$-equivariant sheaf $\calm$, we may decompose it into eigensheaves according to the action of $\mathbb{G}_m$ since $\mathbb{G}_m$ acts trivially on $H$. If $z \in \mathbb{G}_m$ acts on $\calm$ by multiplication by $z^i$ then we say that $\calm$ has weight $i$. The equivariant sheaves of weight $i$ form an abelian subcategory denoted $\qcoh(\caly)_i$. 

We recall the category equivalence between $\qcoh(\caly)_1$ and the right module category $\Mod-\cala'$, where $\cala'$ is any rank $m^2$ Azumaya algebra Morita equivalent to $\cala$. In this case, there is a $GL_m$-equivariant sheaf $\calm$ on $H$ of weight 1 giving $\cala'$ as follows. Note that $\send_H \calm'$ has weight 0 so is just a $PGL_m$-equivariant sheaf. Then $\cala' = (\send_H \calm')^{PGL_m}$ and we may view $\calm'$ as a left $p^*\cala'$-module and a right $\calo_H$-module. The desired category equivalence is given by the inverse functors
$$(-)^{\cala'}:= \shom_H(\calm',-)^{PGL_m}: \qcoh(\caly)_1 \lm \Mod-\cala', \quad
p^*(-)\otimes_{p^*\cala}\ \calm' : \Mod-\cala \lm \qcoh(\caly)_1 .$$

In the special case $\cala' = \cala$, we see from the construction of the $PGL_m$-torsor $p:H \lm Y$ and the framed condition, that $\calm' = \calm \simeq \calo_H^n$. Thus $p^*\cala$ is just the matrix algebra. This gives another relation between $\cala$ and $p$, namely, pulling back by the corresponding $PGL_m$-torsor trivialises the Azumaya algebra to the matrix algebra.

\vspace{2mm}

\section{Endomorphisms compatible with shift}  \label{sendos}  

Let $\calc$ be a $k$-linear category. The two examples of interest for us will be $\Mod-\cala$ where $\cala$ is an Azumaya algebra on a separated algebraic space $Y$ of finite type over $k$, or $\qcoh(\caly)_{1}$ where $\caly = [H/GL_m]$ and $H$ is a $PGL_m$-torsor on $Y$. 

Let $A$ be an $\mathbb{N}$-graded $k$-algebra which is locally finite in the sense that $\dim_k A_i < \infty$ for all $i \in \mathbb{N}$. A {\em (left) $A$-module in $\calc$} is a $\Z$-graded object in $\calc$
$$ \calm_{\bullet} = \oplus_{j \in \Z} \calm_j, \hspace{1cm} \calm_j \in \Obj \calc $$
equipped with multiplication maps $\mu_{ij}: A_i\otimes \calm_j\lm \calm_{i+j}$ which satisfy the usual unit and associativity axioms. Given such an object, there is a natural map $A \lm \End^{\gr}_{\calc} \calm_{\bullet}$ to the graded endomorphism ring $\End^{\gr}_{\calc} \calm_{\bullet} = \oplus_d E_d$ whose $d$-th graded piece is given by 
$$E_d = \prod_{j \in \Z} \Hom_{\calc}(\calm_j,\calm_{j+d}) .$$  
Naturally, $\End^{\gr}_{\calc} \calm_{\bullet}$ is too big to be of much use. I learnt the following trick from James Zhang to reduce to the image of $A \lm \End^{\gr}_{\calc} \calm_{\bullet}$ in a special case.

Suppose now that there is a category auto-equivalence $\s^*: \calc \lm \calc$ and an isomorphism $z:\calm_{\bullet}(1) \stackrel{\sim}{\lm} \s^* \calm_{\bullet}$ of $A$-modules in $\calc$. We write $z_j: \calm_{j+1} \stackrel{\sim}{\lm} \s^* \calm_j$ for the induced isomorphisms. Let $\End_{\calc,z} \calm_{\bullet}$ be the subalgebra of $\End^{\gr}_{\calc} \calm_{\bullet}$ whose $d$-th graded component consists of elements $(\beta_j)_{j \in \Z} \in \prod_{j \in \Z} \Hom_{\calc}(\calm_j,\calm_{j+d})$ such that the following diagram commutes for all $j$.
$$\begin{CD}
\calm_{j+1} @>{\beta_{j+1}}>>  \calm_{d+j+1} \\
@VV{z_j}V  @VV{z_{j+d}}V  \\
\s^*\calm_j @>{\s^*\beta_j}>>  \s^*\calm_{j+d}
\end{CD}$$
This ensures that $\beta_j$ determines $\beta_{j+1}$ and vice versa. Hence we have an isomorphism
$$\End_{\calc,z} \calm_{\bullet} \simeq \oplus_{d \in \Z} \Hom_{\calc} (\calm_0,\calm_d)$$
though on the right hand side, one must be careful to define multiplication by twisting as one finds in the usual twisted homogeneous co-ordinate ring of Artin-Van den Bergh [AV]. Naturally we have an algebra map $A \lm \End_{\calc,z} \calm_{\bullet}$. In the special cases where $\mathcal{C} = \Mod-\cala$ or $\qcoh{\mathcal{S}}$ for some algebraic stack $\mathcal{S}$, we will replace the subscript $\mathcal{C}$ in $\End_{\calc,z}$ and $\End^{gr}_{\calc}$ with $\cala$ or $\mathcal{S}$ accordingly. 

Category auto-equivalences of $\Mod-\cala$ were essentially described by Artin and Zhang [AZ94, corollary~6.9(2)] and involve the notion of an $\cala$-bimodule which we now recall. Let $\calb \in \qcoh(Y \times Y)$ and $p_1,p_2: Y \times Y \lm Y$ be the projection maps. Suppose that the restrictions of $p_1$ and $p_2$ to the support of $\calb$ is finite and that $\calb$ is also given the structure of a $(p_1^* \cala,p_2^*\cala)$-bimodule. Then we say that $\calb$ is a {\em coherent $\cala$-bimodule}. An {\em $\cala$-bimodule} is a direct limit of coherent $\cala$-bimodules. We can tensor with $\calb$ to get a functor as follows $-\otimes_{\cala} \calb:\Mod-\cala \lm \Mod-\cala: \caln \mapsto p_{2*}(p_1^* \caln\otimes_{p_1^* \cala} \calb)$. One can similarly define the tensor product of two $\cala$-bimodules and hence invertible $\cala$-bimodules (see [AZ94, section~6]).

\begin{prop}  \label{pautoequiv}  
Let $\s^*: \Mod-\cala \lm \Mod-\cala$ be an auto-equivalence. Then there is an invertible $\cala$-bimodule $\calb$ such that $\s^*$ is naturally isomorphic to $- \otimes_{\cala} \calb$.
\end{prop}
\textbf{Proof.} The proof for schemes is given in [AZ94, corollary~6.9(2)]. The case of algebraic spaces can be deduced by descent theory.
\vspace{2mm}

We can now calculate $\End_{\calc,z} \calm_{\bullet}$ when $\calc = \qcoh(\caly)_1$. Suppose $\calm_{\bullet}$ is an $A$-module in $\calc = \qcoh(\caly)_{1}$ and that $\calm_0$ has rank $m$. Let $\cala = (\send_H \calm_0)^{PGL_m}$ which we note is Morita equivalent to the Azumaya algebra corresponding the $PGL_m$-torsor $p:H \lm Y$. Suppose further that $\s: \caly \lm \caly$ is a stack automorphism such that $z:\calm_{\bullet}(1) \simeq \s^*\calm_{\bullet}$. We use the category equivalence $(-)^{\cala}: \qcoh(\caly)_{1} \lm \Mod-\cala$ to obtain a category auto-equivalence $\s^*: \Mod-\cala \lm \Mod-\cala$ and an isomorphism $z:\calm_{\bullet}^{\cala}(1) \simeq \s^*\calm_{\bullet}^{\cala}$. Then by proposition~\ref{pautoequiv}, 
there exists an invertible $\cala$-bimodule $\calb$,  such that $\s^*: \Mod-\cala \lm \Mod-\cala$ is given by $- \otimes_{\cala} \calb$. Hence $\calm^{\cala}_d = \calm^{\cala}_0\otimes_{\cala} \calb^{\otimes d} =\calb^{\otimes d}$ and we find
$$ (\End_{\caly,z} \calm_{\bullet})_d = (\End_{\cala,z} \calm_{\bullet}^{\cala})_d = 
\Hom_{\cala}(\calm^{\cala}_0,\calm^{\cala}_d) = H^0(Y,\calb^{\otimes d}).$$
We summarise our discussion in
\begin{prop} \label{ptwistring}  
Let $H \lm Y$ be a $PGL_m$-torsor on a separated algebraic space $Y$ of finite type over $k$, $A$ be a locally finite $\mathbb{N}$-graded algebra and $\calm_{\bullet}$ an $A$-module in $\calc = \qcoh([H/GL_m])_{1}$. Suppose there is a stack automorphism $\s$ of $[H/GL_m]$ such that $\calm_{\bullet}(1) \simeq \s^* \calm_{\bullet}$. Suppose also that $\calm_0$, has minimal rank $m$ and let $\cala$ be the Azumaya algebra $(\send_H \calm_0)^{PGL_m}$. Then there is an invertible $\cala$-bimodule $\calb$ such that the natural map $A \lm \End^{\gr}_{\calc} \calm_{\bullet}$ factors through the twisted ring 
$$ B:= \oplus_{d \in \Z} H^0(Y,\calb^{\otimes d}) .$$
\end{prop}

\vspace{2mm}

\section{Moduli of $m$-points} \label{smodulifat}  

In this section, we construct the moduli stack of simple $m$-points mentioned in the introduction. This stack and the universal simple $m$-point $\calm$ will provide the input data for  proposition~\ref{ptwistring}.

We begin by recalling the basics of the non-commutative projective geometry of Artin-Zhang [AZ94]. Let $A$ be a noetherian $\mathbb{N}$-graded $k$-algebra and $\mmm = A_{>0}$ be the augmentation ideal. We study $A$ via the quotient category $\proj A = A-\Gr/\text{tors}$ where $A-\Gr$ is the category of graded (left) $A$-modules and tors, the localising subcategory of $\mmm$-torsion modules. The relation with projective geometry stems from Serre's theorem: when $A$ is the homogeneous co-ordinate ring of a projective scheme $X$, $\proj A$ is naturally equivalent to the category of quasi-coherent sheaves $\qcoh(X)$. 

We introduce the terminology {\em $A$-tail} to mean an object of $\proj A$. To keep notation unencumbered, we will usually use the same symbol, for example $M$, to denote both an $A$-module and the corresponding $A$-tail though, on rare occasions, we will use $M_{\bullet}$ to denote some module representing the $A$-tail $M$. Similarly, we use $\Ext(-,-)$ for ext groups in $\proj A$ while $\Ext_A(-,-)$ will denote ext groups for $A$-modules. Let $A-\gr, \text{proj}\ A$ denote the full subcategories of $A-\Gr, \proj A$ consisting of noetherian objects. Given a graded module $M\in A-\Gr$ and an integer $n$, we define the shift $M(n)\in A-\Gr$ to be the graded module with $d$-th graded component $M(n)_d := M_{d+n}$. This gives a functor $A-\Gr \stackrel{\sim}{\lm} A-\Gr$ which in turn induces a category isomorphism $\proj A \stackrel{\sim}{\lm} \proj A$. 

Many tools and concepts from projective geometry which can be homologically defined carry over to the non-commutative setting. In particular, we will need to use Artin-Zhang's cohomology functor $H^i(M) = \Ext^i(A,M)$ for $M \in \proj A$. 

We will say a function $h:\mathbb{Z} \lm \mathbb{N}$ has {\em bounded tail} if for some, and hence any $c \in \mathbb{Z}$, the restriction of $h$ to $[c,\infty)$ is bounded.
\begin{defn}  \label{d0dim}  
Let $M$ be an $A$-tail and $M_{\bullet} = \oplus M_n$ a graded module representing it. We define the {\em Hilbert function} of $M$ to be $h:\mathbb{Z} \lm \mathbb{N}: n \mapsto \dim_k H^0(M(n))$ and the {\em Hilbert function} of $M_{\bullet}$ to be $n \mapsto \dim_k M_n$. If $M$ is non-zero and noetherian, then we say $M$ is {\em 0-dimensional} if its Hilbert function  has bounded tail.
\end{defn}
The $0$-dimensional $A$-tails, together with the zero $A$-tail form an abelian subcategory of $\proj A$ which is closed under subquotients. We note that the Hilbert functions of a noetherian $A$-tail and a graded module representing it agree in large degree. 

\subsection{Some geometric assumptions on $A$} \label{ssextraass}  

As in [CN], we will impose some geometric assumptions on our $\mathbb{N}$-graded $k$-algebra $A$. We assume firstly that $A$ is generated in degree one so for any graded $\mmm$-torsionfree $A$-module $M_{\bullet}$, we have $M_0 = 0 \implies M_{-n} = 0$ for $n\in \mathbb{N}$. We further impose the next three conditions on $A$ for the rest of this paper.
\vspace{2mm}

\textbf{1) Strong $\chi$:} $A$ satisfies strong $\chi$ (see [AZ01, (C6.8)] for the rather technical definition). In particular, $A$ is strongly noetherian in the sense that for any commutative noetherian $k$-algebra $R$, $A_R := A \otimes R$ is noetherian. Thus $A$ is noetherian, finitely generated and locally finite. We also know $\proj  A$ is Ext-finite [AZ01, proposition~C6.9] in the sense that for noetherian objects $M \in \proj A, \caln \in \proj A_R$ ($R$ as above), we have $\Ext^i(M,\caln)$ is a noetherian $R$-module. Finally, the strong $\chi$ condition ensures the results of [AZ01, \S E.5] hold. The key result there for us is that the standard Hilbert functor is representable by a scheme called the Hilbert scheme. Furthermore, the Hilbert scheme is a separated, locally finite $k$-scheme which is a countable union of projective schemes [AZ01, theorem~E5.1, E3.1].
\vspace{2mm}

Before describing the next hypothesis, we recall that the local cohomology functor is defined for $M \in A-\gr$ by 
$$H^i_{\mmm}(M) = \lim_{n \lm \infty} \Ext^i_A(A/A_{>n},M).$$

\textbf{2) Finite cohomological dimension:} We assume the cohomological dimension of $A$,
$$ \text{cd}(A) := \sup \{ i | H^i_{\mmm}(-) \neq 0\} $$
is finite. Together with the $\chi$ condition subsumed above, this ensures 
[VdB97] 
the existence of a balanced dualising complex $\w_A \in D^b_c(A\otimes A^{op})$ which gives local duality $H^i_{\mmm}(-)^{\vee} = \Ext^{-i}_A(-,\w_A)$ where $(-)^{\vee}$ is the graded vector space dual. There is also a double Ext spectral sequence namely, for any $M \in D^b_c(A)$, we have 
$$ E^{p,q}_2=\Ext^p_A(\Ext^{-q}_A(M,\w_A),\w_A) \Longrightarrow M .$$
\vspace{2mm}


\textbf{3) CM$_1$:} We assume that for any noetherian graded $A$-module $P_{\bullet}$ whose Hilbert function has bounded tail, we have $H^i_{\mmm}(P_{\bullet}) = 0$ for $i >1$. The notation stands for Cohen-Macaulay in degree one and the hypothesis will hold if $A$ is Cohen-Macaulay with respect to Gelfand-Kirillov dimension in the sense of [AjSZ]. This hypothesis is essentially equivalent to the following weak Grothendieck vanishing result.
\begin{prop} \label{pgrothvan}  
If $P \in \proj A$ is 0-dimensional, then $H^i(P) = 0$ for $i > 0$. 
\end{prop}
\textbf{Proof.} This follows from the fact that if a graded $A$-module $P_{\bullet}$ represents $P \in \proj A$, then $H^i(P) = H^{i+1}_{\mmm}(P_{\bullet})_0$ for $i\geq 1$ by [AZ94, proposition~7.2(2)].
\vspace{2mm}

We also have the following non-vanishing $H^0$ result in dimension 0.
\begin{prop} \label{pnonvanH0}  
If $P \in \proj A$ is 0-dimensional, then $H^0(P) \neq 0$.
\end{prop}
\textbf{Proof.} Let $P_{\bullet}$ be the graded $A$-module $\oplus_{i} H^0(P(i))$ and suppose $H^0(P) = 0$. Now $P_{\bullet}$ is $\mmm$-torsionfree and $A$ generated in degree one so $P_{\bu}$ is generated in positive degrees and hence is noetherian [CN, lemma~4.3]. Our $CM_1$ hypothesis ensures that $H^i_{\mmm}(P_{\bu}) = 0$ for $i > 1$ while [AZ94, proposition~7.2(2)] shows $H^1_{\mmm}(P_{\bu}) = 0$ too. Hence $P_{\bu}$ is $\mmm$-torsion and $P$ must be zero, a contradiction.
\vspace{2mm}

For a commutative noetherian $k$-algebra $R$, we say an $A_R$-tail $\calm \in \proj A_R$, is {\em generated by global sections} if the natural map $H^0(\calm) \otimes A \lm \calm$ is surjective.

\begin{prop}  \label{pgensect}   
Any 0-dimensional $A$-tail $P \in \proj A$ is generated by global sections.
\end{prop}
\textbf{Proof.} We need to show that the map $\phi: H^0(P) \otimes A \lm P$ is surjective so let $C = \coker \phi$. Consider the commutative diagram with exact row
$$\diagram
  &  & H^0(P) \otimes A \dlto \dto  &  &  \\
0 \rto & Q \rto & P \rto & C \rto & 0 
\enddiagram$$
Now $Q$ is 0-dimensional so we must have $H^1(Q) = 0$ and another commutative diagram with exact row
$$\diagram
  &  & H^0(P)\otimes H^0(A) \dlto \dto  &  &  \\
0 \rto & H^0(Q) \rto & H^0(P) \rto & H^0(C) \rto & 0 
\enddiagram$$
Now the vertical map is surjective so $H^0(C) = 0$ and non-vanishing, proposition~\ref{pnonvanH0}, gives $C=0$ as desired.
\vspace{2mm}

\begin{rem}  \label{rCM1}  
{\em We can vary the theory slightly by altering the $CM_1$ hypothesis. Instead of considering graded modules $P_{\bullet}$ whose Hilbert functions have bounded tails, we could have considered a different subcategory of $A-\gr$ which is closed under subquotients. One possibility is those with Gelfand-Kirillov dimension $\leq 1$. This yields a potentially stronger hypothesis, but will also give correspondingly stronger results.}
\end{rem}

We finish this subsection by showing that ``homogeneous co-ordinate rings'' of coherent sheaves of algebras give examples of rings satisfying our three technical hypotheses. To this end, let $Y$ be a projective scheme, $\cala$ be a coherent sheaf of $\oy$-algebras and $\calb$ be an invertible $\cala$-bimodule. We suppose that $\calb$ is {\em ample} in the sense that for any coherent $\cala$-module $\calm$ we have $H^i(Y,\calb^{\otimes n} \otimes_{\cala} \calm) = 0$ for all $i > 0, n \gg 0$. Then one can as usual form the homogeneous co-ordinate ring 
$$A = \bigoplus_{n \geq 0} H^0(Y,\calb^{\otimes n}).$$
\begin{thm}  \label{tordersOK}  
Let $A$ be the homogeneous co-ordinate ring of a coherent sheaf $\cala$ of $\oy$-algebras above. Then $A$ satisfies strong $\chi$, $\text{cd}\, A < \infty$ and $CM_1$. Furthermore, we have $\cala-\Mod \stackrel{\Phi}{\simeq} \proj A$ with the equivalence given by 
$$\calm \mapsto \bigoplus_{n \in \Z} H^0(Y,\calb^{\otimes n}\otimes_{\cala} \calm).$$
\end{thm}
\textbf{Proof.} In the proof of [ASZ, corollary~4.19], they show that in the case $\cala \simeq \oy$, we have $A$ is strongly noetherian and $\qcoh(\oy) \simeq \proj A$. Their proof works verbatim to show the same is true for general $\cala$. In particular, we see that for $\calm \in \cala-\Mod$ we have $H^i(Y,\calm) = \Ext^i_{\cala}(\cala,\calm) = \Ext^i(A,\Phi(\calm))$ so the usual Zariski cohomology for $\cala$-modules equals the Artin-Zhang cohomology of the corresponding $A$-module. Now if $\calm \in \cala-\Mod$ is coherent then we have $\Ext^i_{\cala}(\cala,\calm) = H^i(Y,\calm)$ is finite dimensional for all $i$. Together with ampleness of $\calb$, we see using [AZ94, theorem~7.4] that $A$ satisfies $\chi$. Also, $\Ext^i_{\cala}(\cala,-) = 0$ for $i > \dim Y$. So comparing local and global cohomology using [AZ94, proposition~7.2] we see that $\text{cd}\, A \leq \dim Y + 1$. 

It remains only to show the $CM_1$ condition holds. We first note that a coherent $\cala$-module $\calm$ has finite length if and only if its support as a sheaf on $Y$ is zero-dimensional. For such a module we have the usual (non)-vanishing results, namely, $H^0(Y,\calm) \neq 0$ if $\calm \neq 0$ but $H^i(Y,\calm) = 0$ for $i > 0$. It thus suffices to prove the following 
\begin{lemma} \label{lfinlen}  
With the above notation, the zero-dimensional $A$-tails correspond to the non-zero finite length $\cala$-modules.
\end{lemma}
\textbf{Proof.} We first show that if $\calm \in \cala-\mo$ has finite length then $\dim_k H^0(Y,\calb^{\otimes n} \otimes_{\cala} \calm)$ is bounded as $n$ varies. Now $\calb \otimes_{\cala} -$ is a category auto-equivalence, so it suffices to find a bound on $H^0(Y,\calm)$ as $\calm$ varies over the set of simple $\cala$-modules. This follows since standard generic flatness arguments show that there is a bound on $\dim_k \cala \otimes_Y k(y)$ as $y$ varies over the closed points of $Y$. 

We now prove the converse. First note that by non-vanishing $H^0$ for a simple $\cala$-module $\calm$, we know that the corresponding $A$-module $\Phi(\calm)$ has Hilbert function bounded below by 1. It follows that if $\calm$ is infinite length, then the corresponding $A$-module has unbounded Hilbert function and we are done. 

\subsection{Projectivity of Hilbert schemes of 0-dimensional tails} 
\label{sshilbproj}  

Let $A$ be a graded $k$-algebra generated in degree one satisfying strong $\chi$, $\text{cd} A < \infty$ and $CM_1$ as usual. In preparation for constructing moduli stacks of 0-dimensional $A$-tails, we study in this subsection, auxiliary Hilbert schemes. 

Let $F \in \proj A$ be a noetherian $A$-tail and $\Hilb F$ denote the Hilbert scheme parametrising quotients of $F$. We will call any irreducible closed subscheme of $\Hilb F$ which is maximal with respect to inclusion, an {\em irreducible component}.

\begin{prop}  \label{pcomp}  
The scheme $\Hilb F$ is a countable union of its irreducible components, each of which is projective. Any closed point of $\Hilb F$ lies in finitely many irreducible components. 
\end{prop}
\textbf{Proof.} This follows easily from the fact that $\Hilb F$ is locally finite over $k$ and also a countable union of projective closed subschemes. 
\vspace{2mm}

We first recall Artin-Zhang's notion of flatness [AZ01, section~C.1]. Let $R$ be a commutative noetherian $k$-algebra. For $\calm_{\bullet} \in A_R-\Gr, N \in R-\Mod$, we note that $\calm_{\bullet}\otimes_R N$ is naturally also an object of $A_R-\Gr$. This induces a functor $-\otimes_R -:(\proj A_R)\  \times\ (R-\Mod) \lm \proj A_R$ and we say $\calm \in \proj A_R$ is {\em flat over $R$} if $\calm \otimes_R -$ is exact. A {\em flat family of 0-dimensional $A$-tails over $R$}, is a noetherian object $\calm$ of $\proj A_R$ which is flat over $R$ and for any closed point $\spec K$ of $\spec R$, the induced $A_K$-tail $\calm \otimes_R K$ is 0-dimensional. 

In [AZ01], one partitions the Hilbert scheme using ``tails'' of Hilbert functions. We see below that for 0-dimensional $A$-tails, one can use the Hilbert functions themselves.

\begin{prop} \label{ph0cont}  
Let $\calm$ be a flat family of 0-dimensional $A$-tails parametrised by an affine variety $T$. Then $h^0=\dim_k \circ H^0$ is a constant function on $T$.
\end{prop}
\textbf{Proof.} It was proved in [CN, proposition~5.2ii)] that the Euler characteristic $\chi$ is a continuous function under a regularity assumption on $\proj A$. This assumption is only used to ensure that $h^i=0$ for $i \gg 0$ so $\chi$ is well-defined. However, weak Grothendieck vanishing (proposition~\ref{pgrothvan}) ensures that for any 0-dimensional $A$-tail $M$,  $\chi(M) = h^0(M)$ so the proof in [CN] can also be used to prove the above proposition.
\vspace{2mm}

Given a function $h:\Z \lm \mathbb{N}$ with bounded tail, we let $\Hilb(F,h)$ be the closed subscheme of $\Hilb F$ parametrising quotients $M$ of $F$ with Hilbert function $h$. The above proposition ensures that $\Hilb(F,h)$ is a union of connected components of $\Hilb F$. We now proceed to the main result in this subsection, projectivity of $\Hilb (F,h)$. 

\begin{lemma}  \label{ldualpure}  
Let $M \in A-\gr$ be an $\mmm$-torsionfree module such that $H^i_{\mmm}(M) = 0$ for $i > 1$. Then $H^1_{\mmm}(M)^{\vee}$ is also $\mmm$-torsionfree.
\end{lemma}
\textbf{Proof.} From local duality, we see that the double ext spectral sequence has only one row 
$$ E^{p,1}_2 = \Ext^p_A(\Ext^{-1}_A(M,\w_A),\w_A)  = 
\Ext^p_A(H^1_{\mmm}(M)^{\vee},\w_A).$$
To converge to $M$ we see that 
$$ 0 = \Ext^0_A(H^1_{\mmm}(M)^{\vee},\w_A) = H^0_{\mmm}(H^1_{\mmm}(M)^{\vee})^{\vee}  .$$
This proves the lemma.
\vspace{2mm}

The next theorem strengthens proposition~\ref{pgensect}.
\begin{thm}  \label{tgenSerre}  
Let $M \in \proj A$ be a noetherian $A$-tail of dimension zero. Then the ``Serre module'' $\tilde{M}_{\geq 0} := \oplus_{i\geq 0} H^0(M(i))$ is generated in degree zero. 
\end{thm}
\textbf{Proof.} Let $\tilde{M}_{\bu} := \oplus_{i\in \Z} H^0(M(i))$ and $M_{\bu}$ be the $A$-submodule of $\tilde{M}_{\bu}$ generated by the degree zero part $H^0(M)$. Now the $A$-tail $M$ is generated by global sections by proposition~\ref{pgensect} so the quotient $\tilde{M}_{\bu}/M_{\bu}$ is $\mmm$-torsion. We also know from [AZ94, proposition~7.2(2)] that $\tilde{M}_{\bu}/M_{\bu} = H^1_{\mmm}(M_{\bullet})$. The definition of $M_{\bu}$ ensures also that the degree zero part of the local cohomology group is $H^1_{\mmm}(M_{\bullet})_0 = 0$ so if $M_{\bu} \neq \tilde{M}_{\geq 0}$ then the positive degree part $H^1_{\mmm}(M_{\bullet})_{> 0}$ is a finite dimensional direct summand of $H^1_{\mmm}(M_{\bullet})$. This contradicts the previous lemma that $H^1_{\mmm}(M_{\bullet})^{\vee}$ is $\mmm$-torsionfree.
\vspace{2mm}

\begin{cor}  \label{chilbproj}  
Let $F$ be a noetherian $A$-tail and $h: \Z \lm \mathbb{N}$ be a function with bounded tail. Then the Hilbert scheme $\Hilb(F,h)$ is projective.
\end{cor}
\textbf{Proof.} First note that $F$ can be written as a quotient of $A^m(n)$ for any $m,n \gg 0$. Since $\Hilb(F,h)$ is a closed subscheme of $\Hilb(A^m(n),h)$ we may assume that $F=A^m(n)$ for any $m$ sufficiently large. Shifting $h$ we may also assume $F=A^m$ where $m \geq h(0)$. 

Let $\Hilb_{\gr}$ denote the Hilbert scheme of graded module quotients of $A^m$ with Hilbert function $h$ (in non-negative degrees). We note that $\Hilb_{\gr}$ is projective by [AZ01, theorem~E4.3] and there is a natural morphism $\Hilb_{\gr} \lm \Hilb(A^m,h)$. By proposition~\ref{pcomp}, it suffices to show that the image of $\Hilb_{\gr}$ is dense.

Let $M\in \proj A$ be an $A$-tail with Hilbert function $h$ and $\phi:A^m \lm M$ be a surjection. Let $V$ be the image of $H^0(\phi):H^0(A)^m \lm H^0(M)$ and $T$ be the variety parametrising morphisms $\phi_t: A^m \lm M, t \in T$ such that $\im H^0(\phi_t)$ contains $V$. Also, let $T^0$ the open subset  of $T$ corresponding to morphisms where $H^0(\phi_t)$ is surjective. Our assumption that $m \geq h(0)$ was imposed to ensure that $T^0$ is dense in $T$. Note that there is a flat family of quotients $A^m \otimes \calo_T \lm M \otimes \calo_T$ which when restricted to $t \in T$ is $\phi_t$. This gives rise to a $T$-point of $\Hilb(A^m,h)$ whose image contains the $k$-point corresponding to $\phi$. Now the restricted $T^0$-point lies in the image of $\Hilb_{\gr}$ by theorem~\ref{tgenSerre} so the corollary follows.

\subsection{Cohomology and base change}  \label{sscohbase}  

As usual, in this subsection $A$ will denote a graded $k$-algebra generated in degree one satisfying strong $\chi$, $\text{cd} A < \infty$ and $CM_1$. We study here cohomology and base change in our setting. The proofs mimic the classical ones.

\begin{thm}  \label{tcohbase}   
Let $R$ be a commutative $k$-algebra and $\calm \in \proj A_R$ be noetherian and flat over $R$. Suppose for some $i$, the following condition holds

$(*)$ \hspace{1cm} $H^j(\calm)$ is locally free over $R$ for all $j > i$.

Then there is an isomorphism which is natural in $L \in R-\Mod$
$$(\dag) \hspace{2cm} H^i(\calm) \otimes_R L \simeq H^i(\calm \otimes_R L) \hspace{2cm} .$$
Moreover, if $R$ is noetherian then (*) holds if $H^j(\calm \otimes k(x))=0$ for  all closed points $x \in \spec R$ and $j>i$. In fact, in this case we have $H^j(\calm)=0$.
\end{thm}
\textbf{Proof.} Standard arguments (see for example [Hart, proof of III.12.5]), show that for any right exact functor $F: R-\Mod \lm \mathcal{C}$ to a Grothendieck category $\mathcal{C}$ which preserves direct limits is isomorphic to $F(R) \otimes_R -$. Note that the $F^j:=H^j(\calm \otimes_R -):R-\Mod \lm \proj A$ form a sequence of $\delta$-functors. 

Now our hypothesis that $\text{cd}\ A < \infty$ and [AZ94, proposition~7.2(2)] ensure that $H^d = 0$ for $d \gg 0$. Thus $F^{d-1}$ is right exact and isomorphic to $F^{d-1}(R)\otimes_R -$. If $F^{d-1}(R)$ is locally free then $F^{d-2}$ is right exact too so is isomorphic to $F^{d-2}(R) \otimes_R -$. The isomorphism $(\dag)$ now follows on continuing inductively. 

Suppose now that $R$ is noetherian and $H^j(\calm \otimes k(x))=0$ for  all  $x,j>i$ as above. Ext-finiteness ensures that $H^j(\calm)$ is a noetherian $R$-module so by $(\dag)$ it is zero and the theorem follows.
\vspace{2mm} 

\begin{cor}  \label{cgrflat}  
Let $R$ be a commutative noetherian $k$-algebra and $\calm$ be a flat family of 0-dimensional $A$-tails over $R$. Then $H^i(\calm)=0$ for any $i>0$ and $H^0(\calm)\in R-\mo$ is flat. Furthermore, in this case $H^0$ commutes with base change.
\end{cor}
\textbf{Proof.} Suppose first that $R$ is a finitely generated $k$-algebra. The theorem and weak Grothendieck vanishing proposition~\ref{pgrothvan} ensure that $H^i(\calm) = 0$ for $i>0$. The theorem now shows that $H^0(\calm) \otimes_R - = H^0(\calm \otimes_R -)$ is an exact functor so $H^0(\calm)$ must be flat. The case for general $R$ will follow from the same argument if we can show that for any residue field $K = k(x)$ of a closed point $x \in \spec R$ we have $H^i(\calm \otimes_R K) = 0$ for $i > 0$. To this end we prove 

\begin{lemma}  \label{lintegralform}   
Let $K/k$ be a field extension and $\calm$ be a 0-dimensional $A_K$-tail. There exists a finitely generated $k$-subalgebra $S$ of $K$ and a flat family $\calm_S$ of 0-dimensional $A$-tails over $S$ such that $\calm_S \otimes_S K \simeq \calm$.
\end{lemma}
\textbf{Proof.} Suppose that $\calm$ is represented by the noetherian graded $A_K$-module $M$ and consider a graded free resolution
$$A^r \otimes K \xrightarrow{(x_{ij})} A^s \otimes K \lm M \lm 0 $$
where the matrix co-efficients $x_{ij} \in A_K$. 
There exists a finitely generated $k$-subalgebra $S$ of $K$ such that all the $x_{ij}$ lie in $A_S$. Let $\calm_S$ be the family of $A$-tails over $S$ corresponding to the cokernel of 
$$A^r \otimes S \xrightarrow{(x_{ij})} A^s \otimes S$$
so that $\calm = \calm_S \otimes_S K$. Now $S$ is a domain so generic flatness [AZ01, theorem~C5.1] allows us to replace $S$ with some localisation so that $\calm_S$ satisfies the conditions of the lemma.
\vspace{2mm}

We return now to the proof of the corollary and suppose that $K/k$ is a field extension and $\calm$ is a 0-dimensional $A_K$-tail. Pick a finitely generated subalgebra $S$ of $K$ and flat family $\calm_S$ as in the lemma. Since the corollary has been proved for $S$, we see that for $i>0$ we have $H^i(\calm) = H^i(\calm_S)\otimes_S K = 0$. This completes the proof of the corollary.
\vspace{2mm}

\begin{cor}  \label{crgensect}  
Let $R$ be a commutative noetherian $k$-algebra and $\calm$ be a flat family of 0-dimensional $A$-tails over $R$. Then $\calm$ is generated by global sections. 
\end{cor}
\textbf{Proof.} Suppose first that $R$ is a finitely generated $k$-algebra and consider the exact sequence
$$ H^0(\calm) \otimes A \lm \calm \lm \mathcal{C} \lm 0 $$
where the left hand map is the natural one. Let $x \in \spec R$ be a closed point. Tensoring with $k(x)$ over $R$ gives by corollary~\ref{cgrflat} the exact sequence
$$ H^0(\calm\otimes_R k(x)) \otimes A \lm \calm\otimes_R k(x) \lm \mathcal{C}\otimes_R k(x) \lm 0 .$$
However, $\calm\otimes_R k(x)$ is generated by global sections by proposition~\ref{pgensect} so $\mathcal{C}\otimes_R k(x)$ is zero as must be $\mathcal{C}$ by Artin-Zhang's Nakayama lemma [AZ01, theorem~C3.4].

Again this argument works for general $R$ if we can show that for any field extension $K$ of $k$ then any 0-dimensional $A_K$-tail $\calm$ is generated by global sections. Pick a finitely generated $k$-subalgebra $S\subset K$ such that lemma~\ref{lintegralform} holds so there is a flat family of $A$-tails $\calm_S$ over $S$ with $\calm = \calm_S \otimes_S K$. We have just shown the surjectivity of $H^0(\calm_S) \otimes A \lm \calm_S$. Tensoring with $- \otimes_S K$ and using the fact  that cohomology commutes with base change here by corollary~\ref{cgrflat} shows that $\calm$ is also generated by global sections. 
\vspace{2mm}

For a commutative $k$-algebra $R$, Artin and Zhang give an abstract definition [AZ01, section~B1] of the category $(\proj A)_R$ and later show in proposition~B8.1(v), that it is naturally equivalent to $\proj A_R$. Also, Grothendieck's theory of flat descent works in this setting [AZ01, section~C8], so one can also define $(\proj A)_X$ for $X$ an algebraic space. One sees easily that for a noetherian algebraic space $X$, $(\proj A)_X$ is equivalent to $(A-\Gr)_X/\text{tors}$ where $(A-\Gr)_X$ is the category of $A$-modules in $\qcoh(X)$ and tors is the full subcategory consisting of direct limits of $A$-modules $\calm_{\bullet}$ which are right bounded in the sense that $\calm_n = 0$ for $n \gg 0$. 
 
\subsection{Framed Hilbert schemes}  \label{sshilb}  

We continue our assumption that $A$ is a graded algebra generated in degree one satisfying strong $\chi$, $\text{cd} A < \infty$ and $CM_1$. In this subsection, we examine two open subsets of $\Hilb A^m$ (where $m \in \Z_+$). They will be used to identify the moduli stacks of 0-dimensional $A$-tails. We start by fixing a Hilbert function $h:\mathbb{Z} \lm \mathbb{N}$ with bounded tail and let $m = h(0)$. 

Given an $A$-tail $M$ with Hilbert function $h$, a {\em framing} of $M$ is an isomorphism $k^m \xrightarrow{\sim} H^0(M)$. We can interpret this alternatively as follows. Note first that $k^m$ is naturally a subspace of $H^0(A^m)$. Let $\iota: k^m \hookrightarrow H^0(A^m)$ be the inclusion map. Now linear maps $k \lm H^0(M)$ correspond to morphisms $A \lm M$, so a framing can also be interpreted as a morphism $f:A^m \lm M$ such that $H^0(f)\circ \iota$ is an isomorphism. We hence define the {\em framed part} $\Ht\subset \Hilb(A^m,h)$ of $\Hilb(A^m,h)$ to be the set of quotients $f:A^m \lm M$ where $H^0(f)\circ \iota$ an isomorphism.

\begin{prop}  \label{pframe}  
The framed part $ \Ht$ is an open subset of $\Hilb(A^m,h)$.
\end{prop}
\textbf{Proof.} We need to show that the framed condition is open. Let $R$ be any commutative noetherian $k$-algebra and consider a surjection $\phi: A^m_R \lm \calm$ in $\proj A_R$ where $\calm$ is a flat family of $A$-tails with Hilbert function $h$. 

Define $\mathcal{C}$ to make the following sequence exact.
\[ R^m \lm H^0(\calm) \lm \mathcal{C} \lm 0  \]
where the first map is the natural composite $R^m \lm H^0(A^m_R) \lm H^0(\calm)$. 
For any closed point $x \in \spec R$, we may tensor by $k(x)$ to obtain by corollary~\ref{cgrflat}, the exact sequence
\[ k(x)^m \xrightarrow{f_x} H^0(\calm \otimes_R k(x)) \lm \mathcal{C} \otimes_R k(x) \lm 0 .\]
Hence $f_x$ is a framing precisely when $\mathcal{C} \otimes_R k(x) = 0$, that is, on the complement of the support of $\mathcal{C}$. But $\mathcal{C}$ is a finitely generated $R$-module so we are done.
\vspace{2mm}

Let $R$ be a commutative noetherian $k$-algebra and $\calm$ be a flat family of 0-dimensional $A$-tails over $R$. We say $\calm$ is a family of {\em simple} $A$-tails if for every geometric point $\spec K$ of $\spec R$ we have that $\calm \otimes_R K$ is a simple $A_K$-tail. We 
let $H \subseteq \Ht$ denote the locus corresponding to framed simple $A$-tails with Hilbert function $h$.

\begin{thm}  \label{tsimpleopen}  
Let $R$ be a commutative noetherian $k$-algebra and $\calm$ a flat family of 0-dimensional $A$-tails over $R$. Then the locus of simple $A$-tails in $\spec R$ is open. In particular, $H$ is an open subset of $\Hilb(A^m,h)$.
\end{thm}
\textbf{Proof.} Note corollary~\ref{cgrflat} ensures that $H^0(\calm)$ is a locally free $R$-module of rank $m$. Consider the $\PP^{m-1}$-bundle $\PP_R:= \PP_R(H^0(\calm)^*)$ and the canonical projection morphism $p: \PP_R \lm \spec R$. Let $\call\in \Pic \PP_R$ be the tautological line subbundle of $p^*H^0(\calm)$ and $\Phi$ be the composite map 
\[ \call \otimes_k A \lm p^*H^0(\calm) \otimes_k A \lm p^* \calm  \]
in $(\proj A)_{\PP_R}$. We define $\mathcal{C}$ to be its cokernel. Now Artin-Zhang's version of Nakayama's lemma [AZ01, theorem~C4.3] shows that $\text{Supp}\ \mathcal{C} \subseteq \PP_R$ is closed.  Also, $p$ is proper, so we see it now suffices to prove

\begin{claim}
The locus of $A$-tails which are not simple is $p(\mbox{Supp}\ \mathcal{C})$.
\end{claim}
\textbf{Proof.} For any geometric point $x:\spec K \lm \PP_R$ we have an exact sequence
\[ (\call \otimes A) \otimes_{\PP_R} K \xrightarrow{f_x} p^*\calm \otimes_{\PP_R} K \lm \mathcal{C} \otimes_{\PP_R} K \lm 0  .\]
Consider the induced geometric point $\spec K \xrightarrow{x} \PP_R \xrightarrow{p} \spec R$ of $\spec R$ and the corresponding fibre $\PP_K= \PP_R \times_R K$. Now $x \in \PP_K$ which corresponds to a 1-dimensional subspace $V$ of $H^0(\calm) \otimes_R K = H^0(\calm \otimes_R K)$ by corollary~\ref{cgrflat}. 
Then $f_x$ is essentially just the global section $A_K \lm \calm \otimes_R K$ corresponding to a non-zero element of $V$. 
In particular, $f_x$ is not zero so if $\mathcal{C} \otimes_{\PP_R} K \neq 0$ then $\calm \otimes_R K$ is not simple. Conversely, suppose $\spec K \lm \spec R$ is a geometric point such that $\calm \otimes_R K$ is not simple but has a proper subobject, say $N$. Corollary~\ref{crgensect} ensures that $N$ has a non-zero global section. We consider $K$-points $x:\spec K \lm \PP_K$ of $\PP_K:=\PP_R \times_R K$. We can choose $x$ so it induces a non-zero global section of $N$. For this $x$ we have $\mathcal{C} \otimes_{\PP_R} K  \neq 0$. This completes the proof of the claim and hence the theorem. 
\vspace{2mm}

\subsection{The action of $PGL_m$ on $\Hilb A^m$}  \label{ssglm}  

Let $A$ be as usual a graded algebra generated in degree one satisfying strong $\chi$, $\text{cd} A < \infty$ and $CM_1$. We fix a Hilbert function $h$ with bounded tail and let $m = h(0)$. We know that $GL_m$ acts naturally on $A^m$ and hence on $\Hilb A^m$ and the framed Hilbert scheme $\Ht$ and its simple locus $H$. In this section, we show that $PGL_m$ acts freely on $H$.

We first formulate the action precisely as follows. Note that the pseudo-functor which sends a scheme $X$ to $(\proj A)_X$ yields a category fibred over schemes. Thus given a group $G$ acting on a scheme $X$ we may speak of $G$-equivariant objects in $(\proj A)_X$ (see for example [Vis]). Let $\calm$ denote the universal framed $m$-point on $\Ht$. We construct a natural action of $GL_m$ on $\calm$. Recall that $GL_m$ acts rationally on $k^m$ so there is a coaction $k^m \lm k^m \otimes \calo_{GL_m}$. This induces the first of the maps below, and the second is just the multiplication map on $\calo_{GL_m}$.
$$A^m \otimes \calo_{\Ht} \otimes \calo_{GL_m} \lm A^m  \otimes \calo_{GL_m} \otimes \calo_{\Ht} \otimes \calo_{GL_m} \lm A^m \otimes \calo_{\Ht} \otimes \calo_{GL_m} \lm 
\calm \otimes \calo_{GL_m} = pr^*\calm $$
where $pr: GL_m \times \Ht \lm \Ht$ is the projection. Now the composite above defines a flat family of quotients over $GL_m \times \Ht$ so by the universal property, there is a morphism $\alpha:GL_m \times \Ht \lm \Ht$ which defines the $GL_m$-action on $\Ht$ and satisfies the property that $pr^*\calm \simeq \alpha^*\calm$. This isomorphism defines the $GL_m$-action on $\calm$. We omit the verification of the cocycle condition. Note also that the $GL_m$-action on $\calm$ has weight one. Restricting to $\Gm$, we see that $\Gm$ acts as the identity on $\Ht$ so the action descends to an action of $PGL_m$ on $\Ht$. 

Note that the simple locus $H\subseteq \tilde{H}$ is a $PGL_m$-invariant open subset and the action of $PGL_m$ is set-theoretically free on $H$. We wish to show it is scheme-theoretically free on $H$ under the additional assumption that $H^0(A) = k$. I do not know if this is standard. We will use the following bizarre observation. Let $\s$ be a non-zero $m \times m$-matrix which induces a map $A^m \lm A^m$. If $\bar{H} $ denotes the closure of $H $ in $\Hilb A^m$, then we have an extension of the $GL_m$-action to the following morphism.
$$  (k^{m \times m} - 0) \times H  \lm \bar{H} : (\s, A^m \lm M) \mapsto (\psi:A^m \xrightarrow{\s} A^m \lm M)  .$$
This map is well-defined since if $M$ is simple, it is generated by any of its global sections so the composite $\psi$ on the right is indeed surjective. Note that when $H^0(A) = k$, $\psi$ is not framed unless $\s \in GL_m$. Also as before, changing $\s$ by a scalar has no effect so we obtain a morphism 
$$\alpha: \PP^{m^2-1} \times H  \lm \bar{H}   .$$

\begin{thm}  \label{tfreeact}  
If $H^0(A) = k$, then the group $PGL_m$ acts freely on $H $. In particular, $H/PGL_m$ is a separated algebraic space of finite type.
\end{thm}
\textbf{Proof.} We need to show that the map $(\alpha,p_2): GL_m \times H \lm H \times H$ is a closed immersion, where $p_2$ is the projection. It suffices to show the map is proper, for assuming this, we may apply Zariski's main theorem to see it is finite, and then freeness will follow from the fact that the stabiliser groups of all the closed points are scheme-theoretically trivial. 

To verify the action is proper,  we use the Hilbert-Mumford criterion [MFK, proposition~2.4]. This reduces the verification to 1-parameter subgroups $\chi:\Gm \lm PGL_m$ which we note extends to a regular morphism $\PP^1 \lm \PP^{m^2-1}$. Now $\Hilb_{\gr} A^m$ is constructed in [AZ01, section~E4] as a closed subscheme of a product of Grassmannians, from which it is easy to see the $PGL_m$-action extends to some projective space containing $H$. Hence by [MFK, lemma~0.5], we need only show the map 
$\Phi=(\alpha,p_2): \Gm \times H  \lm H  \times H $
has closed image.

The image is certainly constructible so, arguing by contradiction we assume that we can find a closed curve whose generic point lies in $\im \Phi$ but does not lie globally in $\im \Phi$. More precisely, there exists a smooth projective curve $C$ and a generically 1-1 morphism $f=(f_1,f_2): C \lm \bar{H}  \times \bar{H} $ such that $f(c_0) \in H \times H - \im \Phi$ for some $c_0 \in C$ but $f(c) \in \im \Phi$ for all but finitely many $c \in C$. 

We consider the following commutative diagram whose terms will be defined below.
$$\diagram
B \dto_{\beta} \drto^{\tilde{\phi}} & \\
\PP^1 \times C \drto_{f_2 \circ p_2} \rdashed^{\phi}|>\tip & \bar{H} \times \bar{H} \dto_{p_2} \\
 & \bar{H}
\enddiagram$$
Above, $\phi$ is the composite rational map 
$$\phi = \Phi \circ (1 \times f_2): \PP^1 \times C \lm \PP^1 \times \bar{H} -\,-> \bar{H}  \times \bar{H} $$
which by our bizarre observation above, is a well-defined morphism on $\PP^1 \times f_2^{-1}(H)$.  Resorting to a blowup $\beta$ away from $\PP^1 \times f_2^{-1}(H)$ if necessary, we can resolve the indeterminacy in $\phi$ to obtain the regular morphism $\tilde{\phi}$ above. The image of $\tilde{\phi}$ contains $f(C)$ since the image of $\Phi$ contains the generic point of $f(C)$. Suppose $b \in B$ maps to $f(c_0)$. Note then that if $\beta(b) = (\s,c)$ then $f_2(c) = f_2(c_0) \in H $ so $\beta$ does not blowup any point on $\PP^1 \times c$. Hence $\phi$ restricts to a regular morphism $\phi_c: \PP^1 \times c \lm \bar{H}  \times \bar{H} $  whose image contains $f(c_0)$. 

Now $\chi$ extends to a regular morphism $\tilde{\chi}: \PP^1 \lm \PP^{m^2-1}$. One sees, for example by using the fact that the maximal tori of $PGL_m$ are all conjugate, that for $\s \in \PP^1 - \Gm$ we must have $\tilde{\chi}(\s) \notin PGL_m$. Hence $p_1\circ\phi_c(\s,c)\notin H  $ so $\phi_c(\s,c) \neq f(c_0)$ and there must be some $\tau \in \Gm$ such that $f(c_0) = \phi_c(\tau,c)$. This contradicts the fact that $f(c_0)$ is not in the image of $\Phi$. We have thus shown that $PGL_m$ acts freely.  The last assertion in the theorem follows from [A74, corollary~6.3]. 
\vspace{2mm}

\subsection{Moduli stack of $m$-points}  \label{ssmodstack} 

In this subsection, we construct the moduli stack of 0-dimensional $A$-tails and so finally, obtain our canonical map from $A$ to a twisted ring. As usual, $A$ will be a graded algebra generated in degree one satisfying strong $\chi$, $\text{cd} A < \infty$ and $CM_1$. We also fix a Hilbert function $h:\mathbb{Z} \lm \mathbb{N}$ with bounded tail and let $H$ and $\Ht$ denote the simple framed and framed locus of the Hilbert scheme $\Hilb(A^m,h)$ respectively.

Let $\Aff^{\text{op}}$ denote the category of commutative noetherian $k$-algebras. We will often confuse a pseudo-functor $\caly$ on $\Aff$ with the corresponding pseudo-functor on $\Aff^{\text{op}}$ so for example, $\caly(\spec R) = \caly(R)$.

We wish to define the moduli stack $\tilde{\caly}$ (respectively $\caly$) of (simple) $A$-tails with Hilbert function $h$. For $R \in \Aff^{\text{op}}$ let $\tilde{\caly}(R)$ be the category with objects $\caln \in \proj A_R$ which are flat families of $A$-tails whose Hilbert functions are $h$. In other words, $\caln$ is a noetherian $A$-tail such that $H^0(\caln(i))$ is a locally free $R$-module of rank $h(i)$ for all $i \in \mathbb{Z}$. Similarly, we define $\caly(R)$ to be the category of flat families of simple $A$-tails with Hilbert function $h$. 

Recall from [DM, p.102], [LMB, definition~2.5] that a full subcategory of a stack is a {\em substack} if it is closed under isomorphic objects, closed under pullbacks and the condition of being in the substack is local on $\Aff$ in the flat topology. 

\begin{prop} \label{pHisstack}   
The pseudo-functor $\tilde{\caly}$ defines a stack in groupoids.
\end{prop}
\textbf{Proof.} We know that descent theory for objects in $\proj A$ works fine by [AZ01, section~C8]. We thus a have a stack $\mathcal{S}$ whose category of sections over $R \in \Aff^{\text{op}}$ consists of $R$-flat families of objects in $\proj A$ and with morphisms, the isomorphisms. We wish to show that $\tilde{\caly}$ is a substack of $\mathcal{S}$. 

Now cohomology commutes with base change here by corollary~\ref{cgrflat} so $\tilde{\caly}$ is closed under isomorphisms and pullback. We wish to show that objects descend. Suppose $R \lm R'$ is a faithfully flat map in $\Aff^{\text{op}}$ and let $\caln \in \mathcal{S}(R)$ be such that $\caln \otimes_R R' \in \tilde{\caly}(R)$. Now $H^j(\caln(i)) = 0$ for $j > 0$ and $i \gg 0$ so cohomology of $\caln(i)$ commutes with base change. Thus descent shows $H^0(\caln(i))$ is locally free of rank $h(i)$ so $\caln$ is a flat family of 0-dimensional $A$-tails. Thus cohomology of $\caln$ commutes with base change and $\caln$ must also have Hilbert function $h$.


\vspace{2mm}

\begin{prop}  \label{pHisquot}  
The stack $\tilde{\caly}$ is isomorphic to $[\Ht/GL_m]$ and the open substack $\caly$ is isomorphic to $[H/GL_m]$. 
\end{prop}
\textbf{Proof.} This is fairly standard and we prove the statement about $\tilde{\caly}$ only. First note that there is a natural morphism of stacks 
$$F: \tilde{\caly} \lm [\Ht/GL_m]$$ 
for given $\caln\in \tilde{\caly}(R), R \in \Aff^{\text{op}}$, we can consider the frame bundle $q:\spec \tilde{R} \lm \spec R$ of $H^0(\caln)$ which gives rise to a morphism
\[ \psi:A^m \otimes_k \tilde{R} \lm \caln \otimes_R \tilde{R} .\]
Since $\caln$ is generated by global sections, $\psi$ is surjective and so induces a $GL_m$-equivariant map $\tilde{t}:\spec \tilde{R} \lm  \Ht$. 

We now construct an inverse functor $G: [\Ht/GL_m] \lm \tilde{\caly}$. Consider a cartesian diagram 
$$\begin{CD}
\spec \tilde{R} @>{\tilde{t}}>>   \Ht \\
@VVV @VV{\pi}V  \\
\spec R @>{\tau}>> [ \Ht/GL_m]
  \end{CD}$$
where $\pi$ is the canonical projection. Let $\calm$ be the universal $A$-tail on $\Ht$ with Hilbert function $h$ which we recall is $GL_m$-equivariant. Hence $\tilde{t}^* \calm$ is also $GL_m$-equivariant. But $\tilde{\caly}$ is a stack by the previous proposition so $\tilde{t}^* \calm$ descends to an object $G(\tau) \in \tilde{\caly}(R)$. We conclude that $\tilde{\caly} \simeq [\Ht/GL_m]$.

\vspace{2mm}

Given an $R$-flat family $\caln$ of 0-dimensional $A$-tails with Hilbert function $h(n)$, $\caln(1)$ is also an $R$-flat family of $A$-tails but with Hilbert function $h(n+1)$. If the Hilbert function is constant, then this will lead to a stack automorphism $\s$ of $\tilde{\caly}$ and $\caly$. This motivates the following 

\begin{defn}  \label{dmpoint}  
A noetherian $A$-tail $M$ is called an {\em $m$-point} if $H^0(M(i)) = m$ for all $i \in \Z$. 
\end{defn}

Recall that the natural quotient functor $A-\Gr \lm \proj A$ has a right adjoint $\omega: \proj A \lm A-\Gr$. If $\calm$ is the universal simple $m$-point on $\caly$, then $\calm_{\bullet}:= \omega \calm$ is an $A$-module in $\qcoh(\caly)_1$. We call this the {\em universal simple $m$-point module}. 

The following result is standard.
\begin{lemma}  \label{luniv}  
Let $\s: \caly \xrightarrow{\sim} \caly$ be the stack automorphism given by the shift functor above. The universal simple $m$-point module $\calm_{\bullet}$ satisfies $\s^* \calm_{\bullet} \simeq \calm_{\bullet}(1)$.
\end{lemma}
\textbf{Proof.} Let $\tau: T \lm \caly$ be the morphism defined by a framed simple $m$-point $\caln/T$. We compute $\tau^* \calm \in \qcoh T$ by using the $GL_m$-torsor $\pi: \tilde{T} \lm T$ and $GL_m$-equivariant morphism $\tilde{t}: \tilde{T} \lm H$ which satisfies $\pi^* \caln \simeq \tilde{t}^* \calm$. Then $\tau^* \calm = (\tilde{t}^* \calm)^{GL_m} = \caln$. Hence writing $\calm(\caln/T)$ for $\tau^* \calm$ we see
$$(\s^* \calm)(\caln/T) = \calm(\caln(1)/T) = \caln(1)/T = (\calm(1))(\caln/T). $$
This proves the lemma.

\vspace{2mm}

\begin{thm}  \label{tshift}  
Let $A$ be a $k$-algebra generated in degree one satisfying strong $\chi$, $\text{cd}\, A < \infty, CM_1$ and $H^0(A) = k$. Also, let $\caly=[H/GL_m]$ be the moduli stack of simple $m$-points and $\calm$ be the universal simple on it. If $\cala$ is the Azumaya algebra corresponding to the $PGL_m$-torsor $H \lm Y:=H/PGL_m$, then the natural map $A \lm \End^{gr}_{\caly} \calm_{\bullet}$ factors through the twisted ring 
$$\oplus_{d \in \mathbb{N}} H^0(Y,\calb^{\otimes d})$$
for some invertible $\cala$-bimodule $\calb$.
\end{thm}
\textbf{Proof.} By the lemma, this follows from proposition~\ref{ptwistring} applied to $\calm_{\bullet}$.

\vspace{2mm}
We call the map in the theorem, the {\em canonical map to the twisted ring} induced by the universal simple $m$-point.

\section{Birationally PI} \label{sbirPI}  

We say a semiprime noetherian graded algebra $A$ is {\em birationally PI} if the degree zero part of its graded quotient ring $Q(A)_0$ is PI. By Kaplansky's theorem [Row, theorem~6.1.25], this just means that $Q(A)_0$ is finite over its centre. In this section, we give a criterion for a non-commutative projective surface to be birationally PI in terms of the moduli of $m$-points. Naturally, we hope this will be useful in settling Artin's conjecture on the birational classification of non-commutative projective surfaces. 

We start by looking at an algebraic criterion for being birationally PI. In this section, we let $K$ denote a finite product of fields and $\s$ denote an automorphism of $K$. We consider {\em $\s$-twisted $K$-algebras}, by which we mean a $\Z$-graded algebra $A$ such that $K$ is a subalgebra of $A_0$ and for $\gamma \in K, a_i \in A_i$ we have $a_i \gamma = \s^i(\gamma) a_i$. In Van den Bergh's language of bimodule algebras, this means that $B$ is essentially a $K$-bimodule algebra whose $i$-th graded piece is supported on the graph of $\s^i$. As a first example, let $D$ be a semisimple algebra with centre $K$. Any twisted ring on $D$ (as in proposition~\ref{ptwistring}) is a skew polynomial ring $D[t;\s]$ where $\s$ is an automorphism of $D$ and hence induces an automorphism of $K$. Thus twisted rings on $D$ are $\s$-twisted for appropriate $\s$. Also, subquotients of $\s$-twisted $K$-algebras are $\s$-twisted. 

\begin{prop}  \label{pbiratPI}  
Let $D$ be semisimple algebra which is finite over its centre $K$ and $\s$ be an automorphism of $D$. If $A$ is a semiprime noetherian graded  $k$-subalgebra of $D[t;\s]$ then $A$ is birationally PI. 
\end{prop}
\textbf{Proof.} We let $X$ denote the set of homogeneous regular elements of $A$ so that $Q(A) = AX^{-1}$. If $X$ were also a (right) denominator set in $D[t;\s]$ consisting of regular elements, then we would be done as then $Q(A) \hookrightarrow D[t;\s]X^{-1} = D[t,t^{-1};\s]$. However, this is not clear, so we proceed indirectly via the following lemma.

\begin{lemma}  \label{lOre}  
With notation in the proposition, we have $KA = AK$ which thus is a subalgebra of $D[t;\s]$. Furthermore, $X$ is a (right) denominator set in $KA$. 
\end{lemma}
\textbf{Proof.} The equality $KA = AK$ follows from the fact that $D[t;\s]$ is $\s$-twisted. By Hilbert's basis theorem, $D[t;\s]$ is noetherian so in particular, has ACC on annihilators (of subsets of $D[t;\s]$). It follows that $KA$ also has ACC on annihilators so by [GW, proposition~9.9], it suffices to prove that $X$ is a right Ore set. 

Let $x \in X, \sum \gamma_i a_i \in KA$ where $\gamma_i \in K, a_i \in A$. Since $X$ is an Ore set in $A$ we can find $y_i\in X, b_i \in A$ such that 
$x b_i = a_i y_i$. Now there exists a common right multiple $y \in X$ of the $y_i$'s in $A$ by [GW; lemma~9.2(a)]. Hence, on changing the $b_i$'s, we may assume we have $x b_i = a_i y$ for all $i$. If $\deg x = d$ then 
$$ \sum x \s^{-d}(\gamma_i) b_i = \sum \gamma_i x b_i = (\sum \gamma_i a_i)y$$
which shows that $X$ is an Ore set in $KA$ too.

\vspace{2mm}

We continue the proof of the proposition. The lemma shows that $AX^{-1}$ embeds in $(KA)X^{-1}$ so if $T$ denotes the ideal of $X$-torsion elements in $(KA)$, we see that $A$ embeds in $(KA)/T$. Note also that elements of $X$ remain regular in $(KA)/T$ by definition of $T$ so the proposition will follow from the next lemma applied to $\bar{A} := (KA)/T$. 

\begin{lemma}  \label{lbiratPI}  
Let $\bar{A}$ be a $\s$-twisted $K$-algebra with bounded Hilbert function in the sense that the function $\mathbb{N} \lm \mathbb{N}: n \mapsto \text{length}_K \bar{A}_n$ is bounded. Let $A$ be a semiprime noetherian $\mathbb{N}$-graded $k$-subalgebra such that every homogeneous regular element of $A$ is also regular in $\bar{A}$. Then $A$ is birationally PI.
\end{lemma}
\textbf{Proof.} If $A$ has no regular elements of positive degree then $Q(A)_0 = Q(A_0)$ so the lemma follows since $\bar{A}_0$ (and hence $A_0$) is PI. Let $x\in A$ be an homogeneous regular element of positive degree say $d$. The $K$-bimodule algebra generated by $x$ is of course the skew polynomial ring $K[x;\s^d]$. Now the Hilbert function of $\bar{A}$ is bounded so for $n\gg 0$, multiplication by $x$ gives a surjective map $\bar{A}_n \lm \bar{A}_{n + d}$. Thus $\bar{A}$ is a finitely generated $K[x;\s^d]$-module which must be a direct sum of, say $s$ copies of $K[x;\s]$ possibly shifted (in grading). We see that $\bar{A}$ embeds in the matrix ring $K[x;\s^d]^{s \times s}$.

The lemma will follow if we can show that any homogeneous regular element $y \in A$ is invertible in the ring extension $K[x,x^{-1};\s^d]^{s \times s}$ for then by the universal properties of localisation, $Q(A)$ also embeds in $K[x,x^{-1};\s^d]^{s \times s}$. Let $e_1, \ldots, e_n$ be the primitive idempotents of $K$ so that for each $j$, $Ke_j$ is a field. The sum $\sum_{i \geq 0} K[x;\s^d] e_j y^i$ is not direct so we can find an homogeneous relation of the form 
$$ x^{i_0} \gamma_0 e_jy^0 + x^{i_1} \gamma_1 e_j y + \ldots + x^{i_r} \gamma_r e_j y^r = 0$$
where the $\gamma_l \in K$ and $i_l \in \mathbb{N}$. By assumption, $y$ is also regular in $\bar{A}$ so we may assume that $\gamma_0 e_j\neq 0$ and so has an inverse, say $\gamma'_0$ in $Ke_j$. Then in  $K[x,x^{-1};\s^d]^{s \times s}$ we have $e_j = y'_j y$ where 
$$y'_j = - \gamma'_0 x^{-i_0}(x^{i_1}\gamma_1 e_j + \ldots + x^{i_r} \gamma_r e_j y^{r-1}).$$
Thus $\sum y'_j$ is a left inverse for $y$ in $K[x,x^{-1};\s^d]^{s \times s}$. Also, $K[x,x^{-1};\s^d]^{s \times s}$ has finite Goldie rank so it must be a two-sided inverse by [GW; corollary~5.6]. This completes the proof of the lemma and hence the proposition.

\vspace{2mm}

\begin{lemma}  \label{lsubalgEnd}  
Let $\cala$ be an Azumaya algebra over a separated algebraic space $Y$ of finite type and $\calb$ be an invertible $\cala$-bimodule. Let 
$$ B = \oplus_{d \in \Z} H^0(Y,\calb^{\otimes d}) $$
be the twisted ring as in proposition~\ref{ptwistring}. Then any semiprime noetherian graded subalgebra of $B$ is birationally PI.
\end{lemma}
\textbf{Proof.} First note that $Y$ has an open dense subset which is a scheme by [A74, proposition~4.5]. We may thus restrict $\oplus \calb^{\otimes d}$ to the union $\eta$ of the generic points of $Y$ to obtain a graded algebra say $B_{\eta}$ and a graded algebra homomorphism $B \lm B_{\eta}$ whose kernel is a nilpotent ideal. Thus any semiprime subalgebra of $B$ is a subalgebra of $B_{\eta}$. However, $B_{\eta}$ is a twisted ring on the semisimple algebra $\cala \otimes_Y \calo_{\eta}$ so we are done by proposition~\ref{pbiratPI}.
\vspace{2mm}

 We say that a flat family $\caln/T$ of simple $m$-points is {\em modular} if for every $m$-point $N$, the set $\{ t \in T| \caln \otimes_T k(t)\simeq N\}$ is finite. This is equivalent to saying the map $T \lm Y$ where $Y$ is the coarse moduli space of simple $m$-points is quasi-finite.

\begin{thm}  \label{tbiratPI}  
Let now $A$ be a prime $\mathbb{N}$-graded algebra generated in degree one satisfying strong $\chi, \mbox{cd} A < \infty, CM_1$ and $H^0(A) = k$. Suppose that $gk(A) = 3$ but $A$ has a modular flat family $\caln/T$ of simple $m$-points where $T$ is a 2-dimensional scheme of finite type. Then $A$ is birationally PI. 
\end{thm}
\textbf{Proof.} Let $H$ be the framed Hilbert scheme of simple $m$-points. Since $PGL_m$ acts freely on $H$, we obtain a $PGL_m$-torsor $H \lm Y$ where $Y$ is a separated algebraic space of finite type parametrising simple $m$-points. Let as usual, $\caly = [H/GL_m]$ and $\calm$ be the universal simple $m$-point on $\caly$. From proposition~\ref{ptwistring}, there is an algebra morphism $\Phi:A \lm \End_{\caly,z} \calm$ and the co-domain is a twisted ring. By the lemma above, it suffices to show this map is injective. 

Suppose to the contrary that $I = \ker \Phi \neq 0$. We may replace $T$ with an open subset of an irreducible component of its reduced subscheme and so suppose that $T$ is actually an affine variety. Note that the natural morphism $A \lm \End^{gr}_T \caln$ factors through $\Phi$. Hence $\caln$ is also a modular family of simple $m$-points over $A/I$. Let $P$ be the prime radical of $I$ so that $P/I$ is the semiprime nilpotent ideal in $A/I$. We first show that, on replacing $T$ with some open subscheme if necessary, we can assume $P \caln = 0$ so $\caln$ is even a family of $A/P$-tails. Indeed, by generic flatness [AZ01, theorem~C5.1] we can shrink $T$ until $\caln/P\caln$ is flat. Now for any closed point $t \in T$ we have $P (\caln \otimes_T k(t))=0$ since $\caln \otimes_T k(t)$ is simple and $P$ acts nilpotently. It follows that for all $i \in \Z$ we have a canonical isomorphism $H^0(\caln(i) \otimes_T k(t)) \simeq H^0((\caln/P\caln)(i) \otimes_T k(t))$ so $\caln = \caln /P\caln$ as desired. 

Now $gk(A/P) \leq 2$ so from [AS, theorems~0.3,0.5, corollaries~0.4(ii),9.5], we know that $\proj A/P$ is equivalent to the module category of an order over a projective curve. They do not have 2-dimensional modular families of simple modules so we have obtained our desired contradiction.

\section{A local criterion for being birationally PI} \label{sloccrit}  

In the previous section, we saw that a smooth non-commutative projective surface (with appropriate hypotheses) is birationally PI if it has a surface worth of $m$-points for some $m$. In this section, we use some results of Van den Bergh and Van Gastel [VV] concerning local behaviour of non-commutative surfaces to obtain such a family of $m$-points. Throughout, $A$ will be a prime $\mathbb{N}$-graded algebra generated in degree one satisfying strong $\chi$, $\cd A < \infty, CM_1$ and such that $H^0(A)=k$.

We recall the setup of [VV]. We first assume there is a regular normal homogeneous element $g \in A$ such that $A/gA$ is the twisted ring of a projective Cohen-Macaulay curve $C$. This ensures that the Gelfand-Kirillov dimension of $A$ is three.

Let $p \in C$ be a closed point and $P\in \proj A$ be the corresponding 1-point. We let $\mathcal{C}_p$ denote the full subcategory of $\proj A$ consisting of finite length objects whose composition factors are all isomorphic to $g^i A \otimes_A P$ for various possibly different $i \in \mathbb{Z}$. Using work of Gabriel on locally finite categories, Van den Bergh and Van Gastel show that there is a pseudo-compact ring $\hat{A}_p$ such that $\mathcal{C}_p$ is equivalent to $PCFin(\hat{A}_p)$, the category of finite length pseudo-compact modules in $\hat{A}_p-\Mod$ [VV, theorem~1.1]. We will thus call $\hat{A}_p$ the {\em VV-complete local ring of $A$ at $p$ (or $P$)}. We can take $\hat{A}_p$ to be the endomorphism ring $\End_{\mathcal{C}_p} E$ where $E$ is the direct sum of the distinct injective hulls of the $g^iA \otimes_A P$. The main result of [VV, theorem~1.1] also includes a classification of the possible VV-complete local rings. One possibility is as follows
\begin{equation} \label{eAhatp} 
\hat{A}_p = 
\begin{pmatrix}
R & R & \ldots & R \\
uR & R & & \vdots \\
\vdots & \ddots & \ddots & \\
uR & \hdots & uR & R
\end{pmatrix} \subseteq R^{m \times m}
\end{equation}
where $R = k[[x,y]]$ and $u \in R$ is such that $R/uR\simeq \widehat{\calo}_{C,p}$. In general, $R$ can be a non-commutative regular two-dimensional complete local ring and $m$ is allowed to be infinite. We will not need these facts. Note that there are $m$ simple $\hat{A}_p$-modules and these must correspond to the 1-points $P, gA \otimes_A P, \ldots, g^{m-1} A \otimes_A P$. 

\begin{thm}  \label{tlocallyPI}  
We assume the hypotheses and notation above. Suppose that the VV-complete local ring $\hat{A}_p$ has the form (1) above where $R = k[[x,y]]$. Then $A$ is birationally PI.
\end{thm}
\textbf{Proof.} Let $\tilde{\caly}$ be the stack of $m$-points, $\tilde{H}$ the framed Hilbert scheme of $m$-points and $\caly,H$ the respective simple loci. If we can show that $\dim H \geq m^2+1 = \dim PGL_m + 2$ then all the hypotheses of theorem~\ref{tbiratPI} hold and we may invoke it to conclude $A$ is birationally PI.

The basic idea is very simple. We first consider the $\hat{A}_p \otimes R$-module $R^m$ given by the last column of $\hat{A}_p$ in equation~(1) above. It is flat over $R$ and its simple locus is the Azumaya locus of $\hat{A}_p$, namely $u \neq 0$. This will give a corresponding flat family of $m$-points over $R$ from which we can conclude $\dim H \geq m^2 + 1$. However, since we only have a category equivalence $\mathcal{C}_p \simeq PCFin(\hat{A}_p)$, we must first prove some lemmas. We will need to use several results of Artin-Zhang on adic objects and the reader may wish to refer to [AZ01, section~D2] for the definition and relevant results.
\begin{lemma} \label{ladic}  
Let $S$ be a commutative noetherian complete local $k$-algebra with maximal ideal $\mmm$ such that $S/\mmm^n$ is finite dimensional for all $n$. Then the following categories are all equivalent.
\begin{enumerate}
\item The full subcategory of $\proj A_S$ consisting of noetherian objects $\caln$ such that $\caln \otimes_S S/\mmm^n \in \mathcal{C}_p$ for all $n$.
\item The category of adic objects in $\mathcal{C}_p$.
\item The category of adic objects in $PCFin(\hat{A}_p)$.
\item The full subcategory of $(\hat{A}_p\otimes S)-\mo$ consisting of objects $N$ such that $N \otimes_S S/\mmm^n \in PCFin(\hat{A}_p)$ for all $n$.
\end{enumerate}
\end{lemma}
\textbf{Proof.} The abstract Grothendieck existence theorem of [AZ01, theorem~D6.1] shows that there is an equivalence between the category of noetherian objects in $\proj A_S$ and the category of adic objects in $\proj A$. Restricting this gives the equivalence between i) and ii). Also, ii) and iii) are equivalent since $\mathcal{C}_p \simeq PCFin(\hat{A}_p)$. To see the equivalence between iii) and iv), we first note that here, $\hat{A}_p-\mo$ is equivalent to the category $pc(\hat{A}_p)$ of noetherian pseudo-compact modules. We thus have inverse equivalences given by taking the inverse limit of an adic object and the formal object of a pseudo-compact module. \vspace{2mm}

We use the lemma to show that $R^m$ does indeed correspond to a flat family $\calm$ of $A$-tails over $R$. We check flatness. First note that if $\mmm \triangleleft R$ denotes the maximal ideal, then $R^m \otimes_R R/\mmm^n$ is flat over $R/\mmm^n$ so the same is true of $\calm \otimes_R R/\mmm^n$ by Serre's lemma [AZ01, proposition~C7.1] and the category equivalence $\mathcal{C}_p \simeq PCFin(\hat{A}_p)$. Hence we may invoke [AZ01, theorem~D5.7] to conclude that $\calm$ is also flat over $R$. 
Looking at the closed fibre we see $\calm \otimes_R R/(x,y)$ is an $m$-fold extension of 1-points and hence is an $m$-point. Furthermore, $\calm$ can be represented by a noetherian graded $A \otimes R$-module $\calm_{\bullet} = \oplus \calm_n$ where for $n \gg 0$ the graded pieces $\calm_n$ are free. Thus $\calm$ is a flat family of $m$-points. We obtain thus a morphism $\spec R \lm \tilde{\caly}$. 

We compute now the simple locus.

\begin{lemma} \label{lsimpleM}  
The simple locus of $\calm$ is $u \neq 0$. 
\end{lemma}
\textbf{Proof.} Certainly the simple locus is contained in $u \neq 0$ so we prove the reverse inclusion. Let $x \in \spec R[u^{-1}]$ be a codimension one point of $\spec R$ and $K'$ be the algebraic closure of the residue field $k(x)$. It suffices to show that $\calm \otimes_R K'$ is a simple $(A \otimes K')$-tail so suppose by way of contradiction that $\caln_{K'}$ is a proper subobject. Now $\caln_{K'}$ is noetherian so we can find a finite extension $K/k(x)$ and non-zero $(A\otimes K)$-tail $\caln_K < \calm \otimes_R K$ such that  $\caln_K \otimes_K K' \simeq \caln_{K'}$. Let $S$ be the integral closure of $\im(R \lm K)$ in $K$ so that $S \simeq k[[w]]$. Now corollary~\ref{crgensect} shows that $H^0(\caln_K)$ is a non-zero subspace of $H^0(\calm \otimes_R K) \simeq H^0(\calm \otimes_R S) \otimes_S K$. We may thus choose a global section $A \otimes S \lm \calm \otimes_R S$ whose image $\caln_S$ is such that $\caln_S \otimes_S K$ is a proper subtail of $\calm \otimes_R K$. 

Now the formal object associated to $\calm \otimes_R S$ is an adic object in $\mathcal{C}_p$ so the same is true of the quotient object $\calm \otimes_R S/\caln_S$. Hence we may appeal to lemma~\ref{ladic} to obtain a corresponding quotient $(\hat{A}_p\otimes S)$-module $S^m/N_S$. Now for some $r$, $\hat{A}_p \otimes S$ acts on $S^m$ via the quotient 
$$\hat{A}_{p,S} = 
\begin{pmatrix}
S & S & \ldots & S \\
w^rS & S & & \vdots \\
\vdots & \ddots & \ddots & \\
w^rS & \hdots & w^rS & S
\end{pmatrix}.$$
If $e_1,\ldots ,e_m$ denote the standard idempotents along the diagonal then $N_S = \oplus e_i N_S = \oplus (w^{l_i})$ for appropriate $l_i$. If $l = \max\{l_i\}$ then we see $S^m /N_S$ is an $S/(w^l)$-object. The same must be true of $(\calm\otimes_R S)/\caln_S$ which contradicts the fact that $\caln_S \otimes_S K \neq \calm\otimes_R K$. This concludes the proof of the lemma. \vspace{2mm} 

Let $F\lm \spec R$ be the frame bundle of the locally free module $H^0(\calm)$. We get an induced $GL_m$-equivariant morphism $F \lm \tilde{H}$ as usual. Also lemma~\ref{lsimpleM} ensures the image of $F \times_R R[u^{-1}]$ lies in the open locus $H$ so it suffices to show that the image of $F \times_R R/uR$ has dimension at least $m^2$. Now a direct computation shows that the endomorphism ring $\End \calm \otimes_R R/(x,y)=k$ so the image $F_k$ of $F \times_R R/(x,y)$ in $\tilde{H}$ is isomorphic to $PGL_m$ which has dimension $m^2-1$. It remains to see that the image of $F \times_R R/(u)$ is not (set-theoretically) this copy of $PGL_m$. Now consider a section of $F \times_R R/(u) \lm \spec R/(u)$ and the induced map to $\tilde{H}$. This corresponds to a flat family of quotients of the form $\hat{A}_p^m \otimes R/(u) \lm (R/(u))^m$. Since $(R/(u))^m$ is a non-trivial family over $R/(u)$, we see that a horizontal tangent vector to the fibration $F \lm \spec R$ gets mapped to a tangent vector transverse to $F_k$. This forces the image of $F\times_R R/(u)$ to have dimension at least $m^2$ as desired. The theorem is proved.

\vspace{5mm}
\textbf{\large References}

\begin{itemize}
\item [{[A73]}] M. Artin, ``Th\'eor\`emes de repr\'esentabilit\'e pour les espaces alg\'ebriques'', Les presses de l'universit\'e de Montr\'eal, Montr\'eal, (1973)
  \item [{[A74]}] M. Artin, ``Versal deformations and algebraic stacks'', {\em Invent. Math.} \textbf{27}, (1974), p.165-89
\item [{[A95]}] M. Artin, ``Some problems on three-dimensional graded domains'',  Representation theory and algebraic geometry, Waltham, (1995), London Math. Soc. Lecture Note Ser. \textbf{238}, Cambridge Univ. Press, Cambridge, (1997), p.1-19
\item [{[AS]}] M. Artin, J. T. Stafford, ``Semiprime graded algebras of dimension two''  {\em J. Algebra}  \textbf{227},  (2000),  p.68-123
\item [{[ASZ]}] M. Artin, L. Small, J. Zhang, ``Generic flatness for strongly Noetherian algebras''
 {\em J. Algebra 221} (1999), p.579--610
\item [{[AjSZ]}] K. Ajitabh, S. P. Smith, J. Zhang,``Auslander-Gorenstein rings'' {\em Comm. Algebra}  \textbf{26},  (1998),  p.2159-2180
\item [{[ATV]}] M. Artin, J. Tate, M. Van den Bergh, ``Some algebras associated to automorphisms of elliptic curves'', The Grothendieck Festschrift vol. 1 p.333-85, Prog. Math. \textbf{86}, Birkhauser Boston, (1990)
  \item  [{[AV]}] M. Artin, M. Van den Bergh,  ``Twisted Homogeneous Coordinate Rings'',
       {\em J. of Algebra}, \textbf{133}, (1990) , p.249-271
\item [{[AZ94]}] M. Artin, J. Zhang, ``Noncommutative projective schemes'' {\em Adv. Math.} \textbf{109} (1994), p.228-87
  \item  [{[AZ01]}] M. Artin, J. Zhang,  ``Abstract Hilbert Schemes'', {\em Alg. and Repr. Theory} \textbf{4} (2001), p.305-94
\item [{[CN]}] D. Chan, A. Nyman, ``Non-commutative Mori contractions and $\PP^1$-bundles'', preprint arXiv:0904.1717
\item [{[DM]}] P. Deligne, D. Mumford, ``The irreducibility of the space of curves of a given genus'', {\em Publ. IHES} \textbf{36}, (1969) p.75-110
\item [{[Gom]}] T. Gomez, ``Algebraic stacks''  {\em Proc. Indian Acad. Sci. Math. Sci.}  \textbf{111}  (2001),  p.1-31
\item [{[GW]}] K. Goodearl, R. Warfield, ``An introduction to noncommutative noetherian rings'', London Math. Soc. Texts \textbf{16}, Cambridge University Press, Cambridge, (1989)
\item [{[LMB]}] G. Laumon, Moret-Bailly, ``Champs Alg\'ebriques'', Ergebnisse der Mathematik und ihrer Grenzgebiete. 3. Folge. A Series of Modern Surveys in Mathematics \textbf{39} Springer-Verlag, Berlin, (2000)
\item [{[MFK]}] D. Mumford, J. Fogarty, F. Kirwan, ``Geometric invariant theory'', 3rd enlarged edition, Springer-Verlag, Berlin, (1994)
\item [{[Pat]}] D. Patrick, ``Noncommutative ruled surfaces'' PhD thesis MIT, June 2007
\item [{[Row]}] L. Rowen, ``Ring theory'' student edn, Academic Press, Boston, (1991)
\item[{[RZ]}] D. Rogalski, J. Zhang `` Canonical maps to twisted rings'', {\em Math. Z.}  \textbf{259}  (2008),  p.433-455
\item [{[VdB97]}] M. Van den Bergh, A translation principle for the four-dimensional Sklyanin algebras.  {\em J. Algebra}  \textbf{184},  (1996), p.435-490
  \item  [{[VdB01p]}]   M. Van den Bergh, ``Non-commutative
          $\mathbb{P}^1$-bundles over Commutative Schemes'',
          math.RA/0102005 1 Feb. 2001
\item [{[VV]}] M. Van den Bergh, M. Van Gastel,``Graded modules of Gelfank-Kirillov dimension one of three-dimensional Artin-Schelter regular algebras'', {\em J. Algebra} \textbf{196}, (1997) p.251-82
\item [{[Vis]}] A. Vistoli, ``Grothendieck topologies, fibered categories and descent theory.  Fundamental algebraic geometry'', Math. Surveys Monogr. \textbf{123}, Amer. Math. Soc., Providence, RI, 2005, p.1-104
\end{itemize}

\end{document}